\documentclass[12pt, reqno]{amsart}
\usepackage{amssymb,amsthm,amsmath}
\usepackage[numbers,sort&compress]{natbib}
\usepackage{amssymb,amsmath}
\usepackage{amsfonts}
\usepackage{mathrsfs}
\usepackage{latexsym}
\usepackage{amssymb}
\usepackage{amsthm}
\usepackage{tikz}
\usepackage{tkz-euclide}
\usepackage{float}
\usepackage{pdfsync}
\usepackage{indentfirst}
\date{June 18, 2021}

\usepackage{amsmath}
\usepackage{amsthm}

\renewcommand{\div}{{\rm div}}

\newtheorem{theorem}{Theorem}
\newcommand{\R}{\mathbb R}

\newtheorem{corollary}{Corollary}[section]
\newtheorem{lemma}{Lemma}[section]

\newtheorem{remark}{Remark}[section]

\newcommand{\beq}{\begin{equation}}
\newcommand{\eeq}{\end{equation}}
\newcommand{\ben}{\begin{eqnarray}}
\newcommand{\een}{\end{eqnarray}}
\newcommand{\beno}{\begin{eqnarray*}}
\newcommand{\eeno}{\end{eqnarray*}}

\numberwithin{equation}{section}

\allowdisplaybreaks

\begin{document}
\title[ Asymptotic properties of  the Navier-Stokes equations]{\bf  Asymptotic properties of steady
plane solutions of the Navier-Stokes equations in the exterior of a half-space}
\author{Lili~Wang}
\address[Lili~Wang]{School of Mathematical Sciences, Dalian University of Technology, Dalian, 116024,  China}
\email{wanglili\_@mail.dlut.edu.cn}

\author{Wendong~Wang}
\address[Wendong~Wang]{School of Mathematical Sciences, Dalian University of Technology, Dalian, 116024,  China}
\email{wendong@dlut.edu.cn}
\date{\today}
\maketitle

\begin{abstract}Motivated by Gilbarg-Weinberger's early work on asymptotic properties of steady
plane solutions of the Navier-Stokes equations on a neighborhood of infinity \cite{GW1978} , we investigate asymptotic properties of steady
plane solutions of this system on a half-neighborhood of infinity with finite Dirichlet integral and Navier-slip boundary condition, and obtain that the velocity of the solution grows more slowly than $\sqrt{\log r}$, while the pressure converges to $0$ along each ray passing through the origin.

\end{abstract}



\setcounter{equation}{0}
\section{Introduction}

The
problem of a rigid body through a viscous liquid originates with the pioneering work of Stokes in 1851 \cite{Stokes} on the effect of internal friction
on the movement of a pendulum in a liquid, which is modeled by the incompressible Navier-Stokes equations in a planar exterior
domain (see also the introduction in \cite{Galdi}). When the influence of the boundary walls of the
container was disregarded,
in the fundamental work by Leray \cite{Leray}, where the solution of the two dimensional exterior problem for the Navier-Stokes equations was constructed with finite Dirichlet integral, it is remaining that whether the constructed solution satisfies the asymptotic behavior at $\infty$.
Gilbarg-Weinberger \cite{GW1978} described the asymptotic behavior of the velocity, the pressure and the vorticity, where they showed that $u(x)=o(\ln|x|)$ and
\beno
\lim_{r\rightarrow\infty}\int_0^{2\pi}|u(r,\theta)-\bar{u}|^2d\theta=0
\eeno
for some constant vector $\bar{u}$.
Later, Amick \cite{Amick} proved that $u\in L^\infty$ under zero
boundary condition.  Recently, Korobkov-Pileckas-Russo in \cite{KPR2017} and \cite{KPR2018} obtained that
\beno
\lim_{|x|\rightarrow\infty}u(x)=\bar{u}.
\eeno
More references on the existence and asymptotic behavior of solutions in an exterior domain, we refer to \cite{Ba1973,GNP1997,Ru2009,Ru2010,GG2011,FZ2011,PR2012,KPR2014,DI2017,KR2021} and the references therein. For the generalized q-energy of $\|\nabla u\|_q<\infty$,  asymptotic behavior of solutions or Liouville type theorems are considered in \cite{KTW2021, Wang, Wang2}.
On the other hand, a body moving within an incompressible fluid parallel to a wall in an otherwise unbounded domain is modeled in a planar exterior
domain in a half space with appropriate boundary conditions on the wall. Hillairet-Wittwer \cite{HW2012} obtained the existence and uniqueness of weak solution with the boundary condition at infinity by considering the Oseen model, which implies the faster decay at infinity, see also Boeckle-Wittwer\cite{BW2012}, Guo-Witter-Zhou \cite{GWZ2019} for the recent progress.

In this paper, we investigate the asymptotic properties of Leray solution with finite Dirichlet integral without the infinite boundary condition.
In details,  consider
an arbitrary solution $\left\{\boldsymbol{w}' ,p'\right\}  $ of the following two-dimensional steady Navier-Stokes equations
\ben\label{eq:NS'}	\left\{
\begin{aligned}
\Delta \boldsymbol{w}'-(\boldsymbol{w}'\cdot\nabla)\boldsymbol{w}'&=\nabla p',\\
\nabla\cdot\boldsymbol{w}'&=0\end{aligned}
\right.
\een
in the exterior of the upper half-space, $ \mathbb R_{+}^{2} \setminus \overline{B_{r_0}}$ with $r_0>0$, where $ \boldsymbol{w}'(x,y)=(u'(x,y),v'(x,y)) $ represents fluid velocity, $ p'(x,y) $ is the fluid pressure and $B_{r_0}$ is a ball centered at $0$. 

There are several possibilities for
the boundary condition on $\Gamma=\{(x,y);y=0, |x|>r_0\}$.The widely used is the following Dirichlet or no-slip boundary condition:
\beno
	\boldsymbol{w}'(x,y)=0 \quad \text{on} \quad \partial B_1.  \label{boundary}
\eeno
At this case, Liouville type theorem was proved by Seregin in \cite{Se2015} for the ancient solutions of the time-dependent case. More Liouville type theorems are obtained by Guo-Wang in \cite{GW2022} by assuming different generalized weak solutions. We refer to \cite{SS2011, JSS2013, BS2015, WW2021} and the refereces therein for some results on singularity or classification of solutions under this boundary condition.
However, in the case where the obstacles
have an approximate limit, the Dirichlet boundary conditions are no longer valid (see for example \cite{serrin}). Due to the roughness of the boundary and
the viscosity of the fluid, it is usually assumed that there is a stagnant fluid layer near the boundary, which allows the fluid to slip. This situation seems to match the reality. In 1827, Navier \cite{Navier}
 considered the slip phenomena and proposed the Navier-slip boundary conditions:
\begin{equation}\label{nav boundary}
	\left\{\begin{array}{ll}
		\boldsymbol{w}'(x,y)\cdot n=0,  \\
		2[D(\boldsymbol{w}'(x,y))\cdot n]_{\tau}+\alpha(x)\boldsymbol{w}'_{\tau}=0,
	\end{array}\right.
\end{equation}
where $D(\boldsymbol{w'}(x,y))$ is the stress tensor of fluid,  $n$ and $\tau$ are the unit outer normal vector and tangential vector of the boundary, $\alpha(x)$ is a physical parameter, which
is a $L^\infty $ function on the boundary. For the flat boundary of $\Gamma=\{(x,y);y=0, |x|>r_0\}$, motivated by the above condition of (\ref{nav boundary}), generally  we try to consider
the body suspended in a linear shear flow as in \cite{EBH2020} or \cite{EBH2021}, which showed the forces generated due to wall and shear effects in the absence of relative motion,
\ben\label{eq:slip'}
	v'(x,0)=b , \partial_{2}u'(x,0)=a,\quad {\rm on}~~\Gamma,
\een
where the constants $a,b\in \mathbb{R}$.
Specially, let $\boldsymbol{w}=\boldsymbol{w}'-(ay,b)=(u,v)$ and $p=p'+abx$, then
\ben\label{eq:NS}	\left\{
\begin{aligned}
\Delta \boldsymbol{w}-(\boldsymbol{w}\cdot\nabla)\boldsymbol{w}-(ay\partial_1+b\partial_2)\boldsymbol{w} -(av,0)&=\nabla p,\\
\nabla\cdot\boldsymbol{w}&=0\end{aligned}
\right.
\een
which satisfies the boundary condition:
\ben\label{eq:slip}
	v(x,0)=0 , \partial_{2}u(x,0)=0,\quad {\rm on}~~\Gamma.
\een
Let $\omega=u_{y}-v_{x}$ be the vorticity of a velocity vector $\boldsymbol{w}=(u,v)$ satisfying  (\ref{eq:NS}).
Then the vorticity equation of (\ref{eq:NS})follows:
\begin{equation}\label{eq:vor}
\Delta\omega-\boldsymbol{w}\cdot\nabla\omega-(ay\partial_1+b\partial_2)\omega=0,\quad {\rm in }~~\Omega_0,
\end{equation}
and $\omega=0$ on $\Gamma$ due to (\ref{eq:slip}), where
\beno
\Omega_0=\{(r,\theta), r>r_0,0<\theta<\pi\}=\mathbb{R}_+^2\backslash \overline{B_{r_0}}, \quad r_0>0,
\eeno 
by assuming $r=\sqrt{x^2+y^2}$ with $x=r\cos\theta$ and $y=r\sin\theta$.

Our main results are as follows. The first result is about the growth estimate of the velocity.
\begin{theorem}\label{thm:velocitydecay}
Let $ \boldsymbol{w}\in C^{2}(\overline{\Omega_0}) $ be a solution of the Navier-Stokes equations (\ref{eq:NS}) with the Navier-slip boundary condition (\ref{eq:slip})
and
\ben\label{eq:Dirichlet}
	\int_{\Omega_0}|\nabla\boldsymbol{w}|^{2}dxdy<
	\infty,
\een
then
\begin{equation}\label{eq:velocitydecay}
	\lim_{r\to\infty}\frac{|\boldsymbol{w}(r,\theta)|}{\sqrt{\log r}}=0.
\end{equation}
\end{theorem}
\begin{remark}
The growth rate in this case is similar as that in \cite{GW1978}. However, Cauchy integral formula in \cite{GW1978} doesn't seem to work in this case due to the boundary effect. We used point-wise behavior theorem in \cite{Galdi} by proving the $L^q$ estimate of approximating vorticity functions with $q>2$. It is still unknown whether the convergence of the velocity holds in the $L^\infty$ sense as Korobkov-Pileckas-Russo in \cite{KPR2017} and \cite{KPR2018}. The condition of $ \boldsymbol{w}\in C^{2}(\overline{\Omega_0}) $ can be relaxed to $ \boldsymbol{w}\in C^{2}_{loc}({\Gamma}) $, since the solution with Dirichlet energy (\ref{eq:Dirichlet}) is regular in the interior.
\end{remark}

For the special half-plane, the Liouville type theorem also holds similar to the whole space.
\begin{theorem}\label{thm:liouville}
Let $ \left\{\boldsymbol{w},p\right\} $ be a solution of the Navier-Stokes equations (\ref{eq:NS}) with the Navier-slip boundary condition (\ref{eq:slip}) defined over the upper half space $ \mathbb R_{+}^{2} $. Assume $ \boldsymbol{w}\in C^{2}(\overline{\R_{+}^{2}}) $ and
\ben\label{eq:Dirichletwhole}
\int_{\mathbb R_{+}^{2}}|\nabla\boldsymbol{w}|^{2}dxdy<\infty.
\een
Then $ \boldsymbol{w} $ and $ p $ are constants in $ \mathbb R_{+}^{2}. $
\end{theorem}
\begin{remark}
It's proved by using the monotonicity property of $G(r)=\frac{1}{2}\int_{0}^{\pi}|\boldsymbol{w}(r,\theta)|^{2}d\theta$ and the vanishing vorticity. When the Navier-slip boundary condition (\ref{eq:slip}) is replaced by the no-slip condition, which seems to be more difficult than the high dimensional case in $\mathbb R^{3}$ or $\mathbb R_{+}^{3}$ with the no-slip boundary condition, which is  an open problem as stated in \cite{Galdi} and  \cite{Tsai2018}.  We refer to Galdi \cite{Galdi}, Chae \cite{Chae2014} and Seregin \cite{Seregin}, where
some Liouville type theorems are proved by assuming that $u\in L^{9/2}(%
\mathbb{R}^{3})$, $\Delta u\in L^{6/5}(\mathbb{R}^{3})$ and $u\in BMO^{-1}(\mathbb R^{3})$, respectively.  For more related discussion, we refer to \cite{KPR2015,CW2016,ZhangARMA,WangJDE} and the
references therein.
\end{remark}

Another  major feature of this paper lies in the estimates of the pressure, and we obtained the following results by two different approaches.
\begin{theorem}\label{thm:pressuredecay}
Let $ \left\{\boldsymbol{w},p\right\} $ be a solution of the Navier-Stokes equations (\ref{eq:NS}) with the Navier-slip boundary condition (\ref{eq:slip}). Assume that $ \boldsymbol{w}\in C^{2}(\overline{\Omega_0}),$  (\ref{eq:Dirichlet}) hold, $ a=0 $ and $ bu(r,0)=bu(r,\pi) $ for $ r>r_{0} $.
Then there hold \\
(i)
\begin{equation}\label{eq:pressuredecay2}
	|p(r,\theta)|=o(\log r),
\end{equation}
uniformly for all $\theta\in [0,\pi]$.\\
(ii) up to a constant
\begin{equation}\label{eq:pressuredecay1}
\begin{aligned}
\lim_{r\to\infty,\theta=\theta_{0}}p(r,\theta)=0,
\end{aligned}
\end{equation}
where $\theta_{0}\in [0,\pi].$
\end{theorem}


\begin{remark}
The pressure's decay is not better than the exterior domain of $ \mathbb R^{2}$. The main obstacle comes from the boundary effect. Case (i) of Theorem \ref{thm:pressuredecay} is proved by the help of Brezis-Gallouet inequality, and the proof of Case (ii) in Theorem \ref{thm:pressuredecay} is following the route of \cite{GW1978} but using Green representation formula of harmonic function on the half-annulus domain.
\end{remark}

\begin{remark}
	Combining the two results of Theorem \ref{thm:pressuredecay}, it seems to be unknown that whether  $ p $ converges to zero at infinity along all the rays uniformly. For example, Theorem \ref{thm:pressuredecay} doesn't work when $\theta=\frac{1}{K}$, $ r=K $ and $ p=p(K,\frac{1}{K})=1,$ where $ K\in  \mathbb{N}. $
\end{remark}

At last, we also obtain the decay estimate of the vorticity.
\begin{theorem}\label{thm:vorticitydecay}Let $ \left\{\boldsymbol{w},p\right\} $ be a solution of the Navier-Stokes equations (\ref{eq:NS}) with the Navier-slip boundary condition (\ref{eq:slip}). Assume that $ \boldsymbol{w}\in C^{2}(\overline{\Omega_0})$, (\ref{eq:Dirichlet}) hold and $ a=0 $.
	Then
\begin{equation}\label{eq:voricitydecay}
	\lim\limits_{r\to\infty}\frac{r^{\frac{3}{4}}}{(\log r)^{\frac{1}{8}}}|\omega(r,\theta)|=0,
\end{equation}
where $\theta\in [0,\pi].$
\end{theorem}

The paper is organized as follows, in Section 2, we recall and introduce some technical lemmas. In Section 3, we are aimed to the proof of Theorem \ref{thm:velocitydecay}, which is the decay estimate of the velocity. Theorem {\ref{thm:liouville}}, Theorem \ref{thm:pressuredecay} and Theorem \ref{thm:vorticitydecay} are  proved in Section 4---Section 6, respectively.

Throughout this article, $C$ denotes a constant only depending on the known finite norms of solutions (e.g. $\|\nabla \boldsymbol{w}\|_{L^2}$) and may be different from line to line.

\section{Preliminary lemmas}
The following lemma is similar to the exterior domain $B_{r_0}^c $ of Lemma 2.1 in \cite{GW1978}, one could make the same arguments as \cite{GW1978} in such an exterior angular region. In fact, one may also do an extension to the whole space by noting the following $C^{0,1}$ is enough and we omitted it.
\begin{lemma}\label{lem:meansquare}
Let $ f\in C^{1}(\overline{{\Omega_0}}) $ and have finite Dirichlet integral
\ben\label{eq:Dirichlet-f}
\int_{\Omega_0}|\nabla f|^{2}dxdy<\infty.
\een
Then\ben\label{eq:logconverge}
\lim_{r\to\infty}\frac{1}{\log r}\int_{0}^{\pi}f(r,\theta)^{2}d\theta=0.\een
\end{lemma}

It immediately follows that a special subsequence is pointwise convergent, which is also the same as Lemma 2.2 in \cite{GW1978} and is stated as follows.
\begin{lemma}\label{lem:pointwise}
Let $ \boldsymbol{f}=(f_{1},f_{2}), $ where $ f_{1} $ and $ f_{2} $ satisfy the hypotheses of Lemma \ref{lem:meansquare}. Then there is a sequence $ {r_{n}}, r_{n}\in (2^{n},2^{n+1}), $ such that
\ben\label{eq:pintconverge}
\lim_{n\to\infty}\frac{|\boldsymbol{f}(r_{n},\theta)|^{2}}{\log r_{n}}=0,
\een
where $ \theta\in[0,\pi]. $
\end{lemma}

As in \cite{GW1978}, we can also obtain the bounded-ness of the $L^2$ norm of $\nabla \omega.$
\begin{lemma}\label{lem:Dirichlet-vor}
	 Let $\boldsymbol{w}=(u,v)$ satisfy the Navier-Stokes equations (\ref{eq:NS}) and Navier-slip boundary condition (\ref{eq:slip}) in $\Omega_0$. Assume $ \boldsymbol{w}\in C^{2}(\overline{\Omega_0}) $ and have finite Dirichlet integral of (\ref{eq:Dirichlet}).
	Then
	\begin{equation}\label{eq:Dirichlet-vor}
		\int_{\Omega_0}|\nabla\omega|^{2}dxdy<\infty.
	\end{equation}
\end{lemma}

{\bf Proof.} For a large number $R_0>0$, choose $ r_{0}<r_{1}<\frac{R_0}{2}\leq \rho<  R\leq R_0 $ and a non-negative $ C^{2} $ cut-off function $ \eta(r) $ such that
\begin{equation}\label{eq:cut-off function}
\eta(r)=\left\{\begin{array}{ll}
1,\quad r_{1}<r<\rho,\\0,\quad r\leq r_{0},r\geq R,
\end{array}\right.
\end{equation}
\beno
|\nabla^k\eta|\leq \frac{C}{(r_1-r_0)^k},\quad {\rm as}~~r_0< r< r_1,
\eeno
and
\beno
|\nabla^k\eta|\leq \frac{C}{(R-\rho)^k},\quad {\rm as}~~\rho< r<R,
\eeno
for $k=1,2$.
Let $ h(\omega)=\omega^{2} $ and $ B_{r}^{+}=B_{r}^{+}(0) $ stands for the upper ball centered at $ 0 $ and of radius $ r. $ Then
by $ (\ref{eq:NS})_{2} $ and (\ref{eq:vor}) we have
\begin{equation*}
	\begin{aligned}
	\div&[\eta^6\nabla h(\omega)-h(\omega)\nabla(\eta^6)-\eta^6h(\omega)\boldsymbol{w}]\\&=2\eta^6|\nabla \omega|^{2}-h(\omega)[\Delta(\eta^6)+\boldsymbol{w}\cdot\nabla(\eta^6)]+2\eta^{6}\omega(ay\partial_1+b\partial_2)\omega.
	\end{aligned}
\end{equation*}
 Integration by parts over $ B_{R}^{+}\backslash B_{r_{0}}^{+} $ yields that
\begin{equation}\label{eq:identity}
\begin{aligned}
&\int_{B_{R}^{+}\backslash B_{r_{0}}^{+}}\eta^6|\nabla\omega|^{2}dxdy
\\&=\frac{1}{2}\int_{B_{R}^{+}\backslash B_{r_{0}}^{+}}\omega^{2}(\Delta\eta^6+\boldsymbol{w}\cdot\nabla\eta^6)dxdy+\frac12\int_{B_{R}^{+}\backslash B_{r_{0}}^{+}}\nabla\eta^{6}\cdot(ay,b)\omega^{2}dxdy.
\end{aligned}
\end{equation}
Note that $\eta =1$ for $ r_{1}<r<\rho $ and (\ref{eq:identity}) implies that
\begin{equation}\label{eq:estimation}
	\begin{aligned}
	&\int_{B_{\rho}^{+}\backslash B_{r_{1}}^{+}}|\nabla\omega|^{2}dxdy\\\leq& \int_{B_{R}^{+}\backslash B_{r_{0}}^{+}}|\nabla\omega|^{2}\eta^{6} dxdy\\\leq&\frac{1}{2}\int_{B_{R}^{+}\backslash B_{\rho}^{+}} \omega^{2}(\Delta\eta^{6}+\boldsymbol{w}\cdot\nabla\eta^{6})dxdy+\frac12\int_{B_{R}^{+}\backslash B_{\rho}^{+}}\nabla\eta^{6}\cdot (ay,b)\omega^{2}dxdy\\&+ C(r_1, \|\boldsymbol{w}\|_{L^\infty(B_{r_{1}}^{+}\backslash B_{r_{0}}^{+})},\| \nabla\boldsymbol{w}\|_{L^\infty(B_{r_{1}}^{+}\backslash B_{r_{0}}^{+})} ),
	\end{aligned}
\end{equation}
where we used
\begin{equation}\label{eq:estimation1}
\begin{aligned}
\left|\int_{ {B_{r_{1}}^{+}\backslash B_{r_{0}}^{+}}} \omega^{2}(\Delta\eta^{6}+\boldsymbol{w}\cdot\nabla\eta^{6})dxdy\right|&+\left|\int_{ {B_{r_{1}}^{+}\backslash B_{r_{0}}^{+}}}\nabla\eta^{6}\cdot (ay,b)\omega^{2}dxdy\right|\\&\leq C(r_1, \|\boldsymbol{w}\|_{L^\infty(B_{r_{1}}^{+}\backslash B_{r_{0}}^{+})},\| \nabla\boldsymbol{w}\|_{L^\infty(B_{r_{1}}^{+}\backslash B_{r_{0}}^{+})} ).
\end{aligned}
\end{equation}
First, we consider the second part of the right integral in (\ref{eq:estimation})
\begin{equation}\label{ine:estimate1'}
	\begin{aligned}
	\frac12\int_{B_{R}^{+}\backslash B_{\rho}^{+}}\nabla\eta^{6}\cdot (ay,b)\omega^{2}dxdy\leq \frac{CR}{R-\rho}\int_{B_{R}^{+} \backslash B_{\rho}^{+}}|\nabla\boldsymbol{w}|^{2}dxdy\leq \frac{CR}{R-\rho}.
	\end{aligned}
\end{equation}
The first part of the right integral over the upper half annulus $ B_{R}^{+}\backslash B_{\rho}^{+} $ in (\ref{eq:estimation}) is controlled by
\begin{equation}\label{eq:estimation2}
\begin{aligned}
 &&\left|\int_{B_{R}^{+} \backslash B_{\rho}^{+}}\omega^{2}\Delta\eta^{6}dxdy\right|
 +\left|\int_{B_{R}^{+}\backslash B_{\rho}^{+}}\omega^{2}\left(\boldsymbol{w}-\bar{\boldsymbol{w}}\right)\cdot\nabla\eta^{6} dxdy\right|\nonumber\\&&+\left|\int_{B_{R}^{+}\backslash B_{\rho}^{+}}\omega^{2}\bar{\boldsymbol{w}}\cdot\nabla\eta^{6} dxdy\right|\doteq T_{1}+T_{2}+T_3,
\end{aligned}
\end{equation}
where
$$ \bar{\boldsymbol{w}}(r)=\frac{1}{\pi}\int_{0}^{\pi}\boldsymbol{w}(r,\theta)d\theta. $$
For $ T_{1},$  we get
\begin{equation}\label{eq:estimation3}
\begin{aligned}
T_{1}&\leq\frac{C}{(R-\rho)^{2}}\int_{B_{R}^{+}\backslash B_{\rho}^{+}}\omega^{2}dxdy\leq\frac{C}{(R-\rho)^{2}}\int_{B_{R}^{+}\backslash B_{\rho}^{+}}|\nabla\boldsymbol{w}|^{2}dxdy\leq\frac{C}{(R-\rho)^{2}}.
\end{aligned}
\end{equation}
For $ T_{2}, $ by Schwarz's inequality there holds
\begin{equation*}
	\begin{aligned}
	T_{2}
	&\leq 6\int_{B_{R}^{+}\backslash B_{\rho}^{+}}|\boldsymbol{w}-\bar{\boldsymbol{w}}||\nabla\eta||\eta^{5}|\omega^{2}dxdy
	\\&\leq\frac{C}{R-\rho}\left(\int_{\rho}^{R}\int_{0}^{\pi}|\boldsymbol{w}-\bar{\boldsymbol{w}}|^{2}d\theta rdr\right)^{\frac{1}{2}}\left(\int_{B_{R}^{+}\backslash B_{\rho}^{+}}\eta^{10}\omega^{4}dxdy\right)^{\frac{1}{2}}.
	\end{aligned}
\end{equation*}
 Using Wirtinger's inequality
 $$ \int_{0}^{\pi}|\boldsymbol{w}-\bar{\boldsymbol{w}}|^{2}d\theta\leq \int_{0}^{\pi}|\boldsymbol{w}_{\theta}|^{2}d\theta,$$
one can obtain that
\begin{equation*}
\begin{aligned}
\left(\int_{\rho}^{R}\int_{0}^{\pi}|\boldsymbol{w}-\bar{\boldsymbol{w}}|^{2}d\theta rdr\right)^{\frac{1}{2}}\leq \left(\int_{\rho}^{R}\int_{0}^{\pi}|\boldsymbol{w}_{\theta}|^{2}d\theta rdr\right)^{\frac{1}{2}}\leq R\left(\int_{B_{R}^{+}\backslash B_{r_{0}}^{+}}|\nabla\boldsymbol{w}|^{2}dxdy\right)^{\frac{1}{2}}\leq CR.
\end{aligned}
\end{equation*}
Moreover, since $ \boldsymbol{w}\in C^{2}(\overline{\Omega_0}) ,$ we know that $ \nabla\omega\in C(\overline{\Omega_0}). $  It follows that
$$ \int_{ {B_{r_{1}}^{+}\backslash B_{r_{0}}^{+}}}|\nabla\omega|^{2}dxdy\leq C(\| \nabla^2\boldsymbol{w}\|_{L^\infty(B_{r_{1}}^{+}\backslash B_{r_{0}}^{+})}).$$
 Then by Gagliado-Nirenberg inequality \cite{LL}, we have
\beno
	&&\left(\int_{B_{R}^{+}\backslash B_{\rho}^{+}}\eta^{10}\omega^{4}dxdy\right)^{\frac{1}{2}}
\\&\leq& C\|\eta^{2}\omega\|_{L^{2}(\R_{+}^{2})}\|\nabla(\eta^{2}\omega)\|_{L^{2}(\R_{+}^{2})}\\&\leq& C\|\eta^{2}\omega\|_{L^{2}(\R_{+}^{2})}\|(\nabla\eta^{2})\omega+\eta^{2}(\nabla\omega)\|_{L^{2}(\R_{+}^{2})}\\
&\leq&
	C\|\omega\|_{L^{2}(B_{R}^{+}\backslash B_{r_{0}}^{+})}\left(C+\frac{C}{R-\rho}+\|\nabla\omega\|_{L^{2}(B_{R}^{+}\backslash B_{r_{1}}^{+})}+\|\nabla\omega\|_{L^{2}(B_{r_{1}}^{+}\backslash B_{r_{0}}^{+})}\right)\\
&\leq& C(\| \nabla^2\boldsymbol{w}\|_{L^\infty(B_{r_{1}}^{+}\backslash B_{r_{0}}^{+})})  \left(1+\frac{1}{R-\rho}+\|\nabla\omega\|_{L^{2}(B_{R}^{+}\backslash B_{r_{1}}^{+})} \right).
	\eeno
Hence we obtain the estimate
\ben\label{eq:estimation4}
T_{2}
\leq C(\| \nabla^2\boldsymbol{w}\|_{L^\infty(B_{r_{1}}^{+}\backslash B_{r_{0}}^{+})}) \frac{R}{R-\rho}\left(1+\frac{1}{R-\rho}+(\int_{B_{R}^{+}\backslash B_{r_{1}}^{+}}|\nabla\omega|^{2}dxdy)^{\frac{1}{2}}\right)^{\frac{1}{2}}.
\een
For $ T_{3}, $ it follows from Lemma \ref{lem:meansquare} that $$\bar{\boldsymbol{w}}(r)=o(\sqrt{\log r}),$$ and hence
\begin{equation}\label{eq:estimation5}
\begin{aligned}
T_{3}\leq\int_{B_{R}^{+}\backslash B_{\rho}^{+}}|\bar{\boldsymbol{w}}||\nabla\eta^{6}|\omega^{2}dxdy\leq\frac{C(\log R)^{\frac{1}{2}}}{R-\rho}.
\end{aligned}
\end{equation}
Combining (\ref{eq:estimation})-(\ref{eq:estimation5}), by  Young's inequality  we get
\beno
\int_{B_{\rho}^{+}\backslash B_{r_{1}}^{+}}|\nabla\omega|^{2}dxdy&\leq&\frac12\int_{B_{R}^{+} \backslash B_{r_{1}}^{+}}|\nabla\omega|^{2}dxdy\\&&+C(\| \boldsymbol{w}\|_{C^2(B_{r_{1}}^{+}\backslash B_{r_{0}}^{+})})\left(\frac{R^{2}}{(R-\rho)^{2}}+\frac{(\log R)^{\frac{1}{2}}}{R-\rho}+1\right).
\eeno
Using Giaquinta iteration in \cite{Gia83}, we conclude  that
\beno
\int_{B_{\rho}^{+}\backslash B_{r_{1}}^{+}}|\nabla\omega|^{2}dxdy&\leq& C(\| \boldsymbol{w}\|_{C^2(B_{r_{1}}^{+}\backslash B_{r_{0}}^{+})})\left(\frac{R^{2}}{(R-\rho)^{2}}+1\right).
\eeno
Take $\rho=\frac{R_0}{2}$ and $R=R_0$, and letting $ R_0\to\infty, $ the proof is complete.

\begin{lemma}\label{lem:vorticitydecay}
	 Under the hypotheses of Lemma $ \ref{lem:Dirichlet-vor} $ we have
	\begin{equation}\label{eq:vorticitydecay1}
	\lim_{r\to\infty}r^{\frac{1}{2}}|\omega(r,\theta)|=0,
	\end{equation}
	where $\theta\in[0,\pi].$
\end{lemma}

{\bf Proof.} Using the polar coordinate transformation and Cauchy inequality, we obtain
\begin{equation*}
\begin{aligned}
\int_{2^{n}}^{2^{n+1}}\frac{dr}{r}\int_{0}^{\pi}(r^{2}\omega^{2}+2r|\omega\omega_{\theta}|)d\theta\leq\int_{r>2^{n},0<\theta<\pi}(2\omega^{2}+|\nabla\omega|^{2})dxdy.
\end{aligned}
\end{equation*}
Hence by the integral theorem of the mean, there exists an $ r_{n}\in(2^{n},2^{n+1}) $ such that
\begin{equation}\label{eq:1}
\begin{aligned}
\int_{0}^{\pi}[r_{n}^{2}\omega(r_{n},\theta)^{2}+2r_{n}|\omega(r_{n},\theta)\omega_{\theta}(r_{n},\theta)|]d\theta\leq\frac{1}{\log 2}\int_{r>2^{n},0<\theta<\pi}(2\omega^{2}+|\nabla\omega|^{2})dxdy.
\end{aligned}
\end{equation}
Note that
\begin{equation*}
\begin{aligned}
\omega(r_{n},\theta)^{2}-\frac{1}{\pi}\int_{0}^{\pi}\omega(r_{n},\theta)^{2}d\theta\leq 2\int_{0}^{\pi}|\omega(r_{n},\theta)\omega_{\theta}(r_{n},\theta)|d\theta.
\end{aligned}
\end{equation*}
It follows  from (\ref{eq:1}) and Lemma \ref{lem:Dirichlet-vor} that
\begin{equation*}
\begin{aligned}
0\leq r_{n}\omega(r_{n},\theta)^{2}&\leq\int_{0}^{\pi}[r_{n}^{2}\omega(r_{n},\theta)^{2}+2r_{n}|\omega(r_{n},\theta)\omega_{\theta}(r_{n},\theta)|]d\theta\\&\leq\frac{1}{\log 2}\int_{r>2^{n},0<\theta<\pi}(2\omega^{2}+|\nabla\omega|^{2})dxdy\to 0,  n\to\infty,
\end{aligned}
\end{equation*}
which implies
\begin{equation}\label{eq:max}
\begin{aligned}
\lim_{n\to\infty}[r_{n}\max\limits_{\theta\in[0,\pi]}\omega(r_{n},\theta)^{2}]=0.
\end{aligned}
\end{equation}
By (\ref{eq:slip}), we have $ \omega(x,0)=0 $ on $\Gamma.$ Let $ D $ denote the half-annulus region with $r_n<r<r_{n+1}$.  Since $\omega$  a solution of the equation (\ref{eq:vor}),  it satisfies the maximum principle in $ D. $ Noting that $ r_{n+1}\leq 4r_{n}, $ we infer that for $ r\in(r_{n},r_{n+1}) $
\begin{equation*}
\begin{aligned}
r\max_{\theta\in [0,\pi]}\omega (r,\theta)^{2}
&\leq r_{n+1}\max_{\theta\in [0,\pi]}\omega (r,\theta)^{2}\\&\leq r_{n+1}\max[0,\max_{\theta\in [0,\pi]} \omega(r_{n},\theta)^{2},\max_{\theta\in [0,\pi]}\omega(r_{n+1},\theta)^{2}]\\&
\leq\max[0,4r_{n}\max_{\theta\in [0,\pi]}\omega(r_{n},\theta)^{2},r_{n+1}\max_{\theta\in [0,\pi]}\omega(r_{n+1},\theta)^{2}],
\end{aligned}
\end{equation*}
which implies the desired result (\ref{eq:vorticitydecay1}) due to (\ref{eq:max}).

At last, we introduce the Brezis-Gallouet inequality (see Lemma 2 in \cite{BG1980}, or Lemma 3.1 in \cite{CPZ2020}).
\begin{lemma}\label{lem:BGI}
	Let $ f\in H^{2}(\Omega), $ where $\Omega$ is a bounded domain or an exterior domain with compact smooth boundary. Then there exists a constant $ C_{\Omega} $ depending only on $\Omega$, such that
	$$\left\|f\right\|_{L^{\infty}(\Omega)}\leq C_{\Omega}\left\|f\right\|_{H^{1}(\Omega)}\log^{\frac{1}{2}}(e+\frac{\left\|\Delta f\right\|_{L^{2}(\Omega)}}{\left\|f\right\|_{H^{1}(\Omega)}}),$$
	or
	$$\left\|f\right\|_{L^{\infty}(\Omega)}\leq C_{\Omega}(1+\left\|f\right\|_{H^{1}(\Omega)})\log ^{\frac{1}{2}}(e+\left\|\Delta f\right\|_{L^{2}(\Omega)}).$$
\end{lemma}

Note that the second inequality can be obtained immediately from the first one by arguments whether $ \left\|f\right\|_{H^{1}(\Omega)}<1 $.

The following lemma is from Theorem II.9.1 in \cite{Galdi}.
\begin{lemma}\label{lem:convergence-f}
	Let $\Omega \subset\R^{2}$ be an exterior domain and let
	$$\nabla f\in L^{2}\cap L^{p}(\Omega),$$ for some $ 2<p<\infty. $ Then
	$$\lim_{|x|\to\infty}\frac{|f(x)|}{\sqrt{\log (|x|)}}=0,$$
	uniformly.
\end{lemma}

\section{Deacy of the velocity}

{\bf Proof of Theorem \ref{thm:velocitydecay}.} For  $ r=\sqrt{x^{2}+y^{2}}, $ we take two cut-off functions $\varphi(x,y)$ and $\eta(x,y)$ as follows:
\begin{equation*}
	\varphi(r)=\left\{
	\begin{aligned}
	&1,\quad r<R\\&0,\quad r>2R
	\end{aligned}
	\right.,\quad
	\eta(r)=\left\{
	\begin{aligned}
	1,\quad r>3r_{0}\\0,\quad r<2r_{0}
	\end{aligned}.
	\right.
\end{equation*}

Firstly, we show that
\begin{equation}\label{eq:vor-bound}
	\omega\in L^{p}(\Omega_0),\quad \forall p\geq 2.
\end{equation}
Noting that $|\nabla\eta|\leq C$ and $ |\nabla\varphi|\leq\frac{C}{R} $ for a constant $ C $ independent of $ R. $
Then using (\ref{eq:Dirichlet-vor}), (\ref{eq:Dirichlet}) and Gagliardo-Nirenberg inequality, for any $p\geq 2,$ we have
\begin{equation*}
	\begin{aligned}
	\left\|\omega\varphi\eta\right\|_{L^{p}(\R_{+}^{2})}&\leq C\left\|\omega\varphi\eta\right\|_{L^{2}(\R_{+}^{2})}^{\frac{2}{p}}\left\|\nabla(\omega\varphi\eta)\right\|_{L^{2}(\R_{+}^{2})}^{1-\frac{2}{p}}\\&\leq C\left\|\omega\right\|_{L^{2}(\R_{+}^{2}\backslash B_{2r_{0}}^{+})}^{\frac{2}{p}}\left(\left\|(\nabla\omega)\varphi\eta\right\|_{L^{2}(\R_{+}^{2})}^{1-\frac{2}{p}}+\left\|(\nabla\varphi)\omega\eta\right\|_{L^{2}(\R_{+}^{2})}^{1-\frac{2}{p}}+\left\|(\nabla\eta)\omega\varphi\right\|_{L^{2}(\R_{+}^{2})}^{1-\frac{2}{p}} \right)\\&\leq C\left\|\nabla\boldsymbol{w}\right\|_{L^{2}(\R_{+}^{2}\backslash B_{2r_{0}}^{+})}^{\frac{2}{p}}\left(C\left\|\nabla\omega\right\|_{L^{2}(\R_{+}^{2}\backslash B_{2r_{0}}^{+})}^{1-\frac{2}{p}}+\frac{C}{R}+C\right)\\&\leq C(1+R^{-1}).
	\end{aligned}
\end{equation*}
Since $\varphi=1$ as $ R\to\infty. $ Hence let $ R\to\infty, $ we get
$$\left\|\omega\right\|_{L^{p}(\R_{+}^{2}\backslash B_{3r_{0}}^{+})}\leq C.$$
Moreover, due to $ \boldsymbol{w}\in C^{2}(\overline{\Omega_0}), $ there holds $\omega\in C^{1}(\overline{\Omega_0}).$ Then there holds
\begin{equation*}
	\left\|\omega\right\|_{L^{p}(B_{3r_{0}}^{+}\backslash B_{r_{0}}^{+})}\leq C.
\end{equation*}
Thence (\ref{eq:vor-bound}) holds.

Secondly, set $ \hat{\boldsymbol{w}}(x,y)=(\hat{u}(x,y),\hat{v}(x,y)), $ where
\begin{equation*}
	\hat{u}(x,y)=\left\{\begin{aligned}
	u(x,y),~&y\geq 0\\u(x,-y),~&y<0
	\end{aligned}
	\right.,~~\hat{v}(x,y)=\left\{\begin{aligned}
	v(x,y),~&y\geq 0\\-v(x,-y),~&y<0
	\end{aligned}
	\right..
\end{equation*}
Then
\begin{equation*}
	\hat{\boldsymbol{w}}\in C_{loc}^{0,1}(\R^{2}\backslash B_{r_{0}}),
\end{equation*}
and
\begin{equation*}
	\hat{\boldsymbol{w}}=\boldsymbol{w},\quad {\rm~in~} \Omega_0.
\end{equation*}
Moreover
\begin{equation*}
	\hat{\omega}(x,y)=\left\{\begin{aligned}
\partial_{y}u(x,y)-\partial_{x}v(x,y),~y\geq 0\\-(\partial_{y}u)(x,-y)+(\partial_{x}v)(x,-y),~y<0
	\end{aligned}
	\right..
\end{equation*}
Let $\varphi$ is standard mollifier operator, $\varepsilon>0,$ set $ \varphi_{\varepsilon}=\varepsilon^{-2}\varphi(\frac{x}{\varepsilon}). $ Thus
$$ \varphi_{\varepsilon}*\hat{\boldsymbol{w}}=\hat{\boldsymbol{w}_\varepsilon}\in C^{\infty}(\R^{2}\backslash B_{2r_{0}}).$$
For the cut-off function $\phi = 1$ outside $B_{3r_0}$ and $\phi = 0$ in $B_{2r_0}$,
noting that
\begin{equation*}
	\nabla^\bot \cdot (\hat{\boldsymbol{w}_\varepsilon} \phi) =\hat{\omega_\varepsilon} \phi + \hat{\boldsymbol{w}_\varepsilon}\cdot \nabla^\bot \phi,
\end{equation*}
where
\begin{equation*}
	\hat{\omega_\varepsilon} = \nabla^\bot \cdot \hat{\boldsymbol{w}_\varepsilon}.
\end{equation*}
Putting the operation $\nabla^\bot$ into the both side of the above equation, we have
\begin{equation}\label{eq:extension-vel1}
	\Delta (\hat{\boldsymbol{w}_\varepsilon} \phi) = \nabla^{\bot} (\hat{\omega_\varepsilon} \phi + \hat{\boldsymbol{w}_\varepsilon} \cdot \nabla^\bot \phi)+\nabla[\nabla\cdot(\hat{\boldsymbol{w}_{\varepsilon}}\phi)].
\end{equation}
Noting that
\begin{equation}\label{eq:extension-vel2}
	\nabla (\hat{\boldsymbol{w}_\varepsilon}\phi) = \nabla \hat{\boldsymbol{w}_\varepsilon} \to 0 \quad {\rm as} \quad r \to \infty,
\end{equation}
by Calderon-Zygmund estimates, we have
\begin{equation*}
	\|\nabla (\hat{\boldsymbol{w}_\varepsilon} \phi)\|_{L^p(\mathbb{R}^2)} \leq C \|(\hat{\omega_\varepsilon} \phi + \hat{\boldsymbol{w}_\varepsilon}\cdot \nabla^\bot \phi)\|_{L^p(\mathbb{R}^2)}+C\|\nabla\cdot(\hat{\boldsymbol{w}_{\varepsilon}}\phi)\|_{L^{p}(\mathbb{R}^{2})},\quad p>1.
\end{equation*}
Noting that $ \nabla\cdot\hat{\boldsymbol{w}_{\varepsilon}}=0, $ then by triangle inequality, we have
\begin{equation*}
	\|\nabla (\hat{\boldsymbol{w}_\varepsilon}\phi)|_{L^p(\mathbb{R}^2)} \leq C \|\hat{\boldsymbol{w}_\varepsilon}\|_{L^p(B_{3r_0} \setminus B_{2r_0})} + C\|\hat{\omega_\varepsilon} \phi\|_{L^p(\mathbb{R}^2)}.
\end{equation*}
Since $\hat{\boldsymbol{w}_\varepsilon} \in C^{0,1}(\R^{2}\backslash B_{2r_{0}})$ uniformly, there holds  $\|\hat{\boldsymbol{w}_\varepsilon}\|_{L^p(B_{3r_0} \setminus B_{2r_0})} \leq C$. Using the definition of $\phi$, we have
\begin{equation*}
	\|\nabla \hat{\boldsymbol{w}_\varepsilon}\|_{L^p(\mathbb{R}^2 \setminus B_{3r_0})} \leq C + C \|\hat{\omega_\varepsilon} \phi\|_{L^p(\mathbb{R}^2)} \leq C + \|\hat{\omega} \|_{L^p(\mathbb{R}^2 \setminus B_{2r_0})}.
\end{equation*}
Let $\varepsilon \to 0$, by Lebesgue's dominated convergence theorem and (\ref{eq:vor-bound}), we have
\begin{equation*}
\begin{aligned}
\|\nabla \hat{\boldsymbol{w}}\|_{L^p(\mathbb{R}^2 \setminus B_{3r_0})} &\leq C + C \|\hat{\omega} \|_{L^p(\mathbb{R}^2 \setminus B_{2r_0})}\\&\leq C+C\left\|\omega\right\|_{L^{p}(\R_{+}^{2}\backslash B_{2r_{0}}^{+})}\leq C,\quad \forall~~ p\geq 2.
\end{aligned}
\end{equation*}
Thus with the help of Lemma \ref{lem:convergence-f}, we have
\begin{equation*}
	\hat{\boldsymbol{w}}=o(\sqrt{\log r}),
\end{equation*}
which implies
\begin{equation*}
	\boldsymbol{w}=o(\sqrt{\log r}).
\end{equation*}
The proof is complete.

\section{Liouville type  theorem in the half-plane}

{\bf Proof of Theorem {\ref{thm:liouville}}.} In the case of $ \mathbb R_{+}^{2}$, it follows from Lemma \ref{lem:vorticitydecay} that $\omega\to 0$ at infinity. By (\ref{eq:slip}), $\omega(x,0)=0.$ Hence $\omega\equiv 0$ in $ \mathbb R_{+}^{2} $ due to the maximum principle. Noting that $ \Delta u=\omega_{y} $ and $ \Delta v=-\omega_{x}, $ there holds  $ \Delta\boldsymbol{w}=(\Delta u,\Delta v)=0 $ in $ \mathbb R_{+}^{2}. $

Recall $ B_{r}^{+}=\left\{\sqrt{x^{2}+y^{2}}<r : y> 0\right\} $ and $ \partial B_{r}^{+}=\Gamma_{1}+\Gamma_{2}, $ where \beno \Gamma_{1} =\left\{(x,y):-r\leq x\leq r,y=0\right\},\quad  \Gamma_{2}=\left\{(x,y):\sqrt{x^{2}+y^{2}}=r,y>0\right\} . \eeno
 Due to  (\ref{eq:slip}), we get
\begin{equation*}
	\int_{\Gamma_{1}}\boldsymbol{n}\cdot\nabla\boldsymbol{w}\cdot\boldsymbol{w}ds=-\int_{-r}^{r}\frac{\partial}{\partial y}u(x,0)\cdot u(x,0)dx=0.
\end{equation*}
Then by the integration by parts formula, we have
\begin{equation}\label{eq:liouville-1}
\begin{aligned}
0&=\int_{B_{r}^{+}}-\Delta\boldsymbol{w}\cdot\boldsymbol{w}dxdy\\&=\int_{B_{r}^{+}}|\nabla\boldsymbol{w}|^{2}dxdy-\int_{\Gamma_{1}+\Gamma_{2}}\boldsymbol{n}\cdot\nabla\boldsymbol{w}\cdot\boldsymbol{w}ds\\&=\int_{B_{r}^{+}}|\nabla\boldsymbol{w}|^{2}dxdy-\frac{1}{2}\int_{0}^{\pi}\frac{\partial}{\partial r}[|\boldsymbol{w}(r,\theta)|^{2}]d\theta\cdot r,
\end{aligned}
\end{equation}
where $ \boldsymbol{n} $ denotes outward unit normal to $ \partial B_{r}^{+}. $
Write
$$ G(r)=\frac{1}{2}\int_{0}^{\pi}|\boldsymbol{w}(r,\theta)|^{2}d\theta,$$
and
$$ C_{0}=\int_{B_{r_{0}}^{+}}|\nabla \boldsymbol{w}|^{2}dxdy,$$ where $ r>r_{0}. $ Then
$$ rG^{'}(r)
\geq C_{0}$$ for any $ r>r_{0}. $
Integration it over $ (r_{0},r), $ we have
$$ G(r)\geq G(r_{0})+C_{0}\log\frac{r}{r_{0}}.$$
Lemma \ref{lem:meansquare} shows that $ G(r)=o(\log r), $ hence
$$ C_{0}=\int_{B_{r_{0}}^{+}}|\nabla\boldsymbol{w}|^{2}dxdy=0.$$
Letting $ r_{0}\to\infty, $ it follows that $\nabla \boldsymbol{w}=\boldsymbol{0} $ in $ \mathbb R_{+}^{2} $ and  $ \boldsymbol{w} $ is constant.The proof is complete.

\section{Asymptotic behavior of the pressure}



Since $ \Delta u=\omega_{y} $ and $ \Delta v=-\omega_{x}, $ the Navier-Stokes equations can be written in the form
\begin{equation}
\left\{\begin{aligned}
\omega_{y}+(u+ay)v_{y}-(v+b)u_{y}-av&=p_{x},\\-\omega_{x}-(u+ay)v_{x}+(v+b)u_{x}&=p_{y},\\u_{x}+v_{y}&=0.
\end{aligned}\right.
\end{equation}
Specially, let $ a=0, $ we get
\begin{equation}\label{eq:component-NS}
\left\{\begin{aligned}
\omega_{y}+uv_{y}-(v+b)u_{y}&=p_{x},\\-\omega_{x}-uv_{x}+(v+b)u_{x}&=p_{y},\\u_{x}+v_{y}&=0,
\end{aligned}\right.
\end{equation}
it follows that
\begin{equation}\label{eq:partial pressure}
\begin{aligned}
p_{r}=\frac{1}{r}[\omega_{\theta}+uv_{\theta}-(v+b)u_{\theta}].
\end{aligned}
\end{equation}

The first lemma is about the convergence of the square norm for a subsequence, which is simialr as in \cite{GW1978}.
\begin{lemma}\label{lem:pressure convergence1}
	Let $ \left\{\boldsymbol{w},p\right\} $ be a solution of the Navier-Stokes equations (\ref{eq:NS}) with the Navier-slip boundary condition (\ref{eq:slip}). Assume that $ \boldsymbol{w}\in C^{2}(\overline{\Omega_0}),$  (\ref{eq:Dirichlet}) hold and $ a=0, $
	 there exists a sequence $ \left\{R_{n} \right\}, R_{n}\in (2^{2^{n}},2^{2^{n+1}}), $ such that
	\begin{equation}\label{eq:pressure convergence1}
	\begin{aligned}
	\lim_{n\to\infty}\int_{0}^{\pi}|p(R_{n},\theta)-\bar{p}(R_{n})|^{2}d\theta=0.
	\end{aligned}
	\end{equation}
\end{lemma}

 {\bf Proof.} We first show that for any $ r_{1}>\max (r_{0},1) $
$$ \int_{r>r_{1},0<\theta<\pi}\frac{|\nabla p|^{2}}{\log r}dxdy<\infty.$$
Using $ (\ref{eq:component-NS})_{1,2} $ and Cauchy inequality
\begin{equation*}
	\begin{aligned}
	|\nabla p|^{2}&=p_{x}^{2}+p_{y}^{2}\\&=[\omega_{y}+uv_{y}-(v+b)u_{y}]^{2}+[-\omega_{x}-uv_{x}+(v+b)u_{x}]^{2}\\&\leq 4|\nabla\omega|^{2}+8u^{2}|\nabla v|^{2}+4(2v^{2}+2b^{2})|\nabla u|^{2}\\&\leq 4|\nabla\omega|^{2}+16|\boldsymbol{w}|^{2}|\nabla\boldsymbol{w}|^{2}+8b^{2}|\nabla\boldsymbol{w}|^{2}.
	\end{aligned}
\end{equation*}
From (\ref{eq:velocitydecay}) in Theorem \ref{thm:velocitydecay}, we know $ |\boldsymbol{w}(r,\theta)|^{2}=o({\log r}) $ , thus
\begin{equation*}
\begin{aligned}
&\int_{r>r_{1},0<\theta<\pi}\frac{|\nabla p|^{2}}{\log r}dxdy\\&\leq 4 \int_{r>r_{1},0<\theta<\pi}\frac{|\nabla\omega|^{2}}{\log r}dxdy+16\int_{r>r_{1},0<\theta<\pi}\frac{|\boldsymbol{w}|^{2}}{\log r}|\nabla\boldsymbol{w}|^{2} dxdy+8b^{2}\int_{r>r_{1},0<\theta<\pi}\frac{|\nabla\boldsymbol{w}|^{2}}{\log r}dxdy\\&\leq \frac{4}{\log r_{1}}\int_{r>r_{1},0<\theta<\pi}|\nabla\omega|^{2}dxdy+C\int_{r>r_{1},0<\theta<\pi}|\nabla\boldsymbol{w}|^{2}dxdy+\frac{8b^{2}}{\log r_{1}}\int_{r>r_{1},0<\theta<\pi}|\nabla\boldsymbol{w}|^{2}dxdy\\&\leq C\left(\int_{r>r_{1},0<\theta<\pi}|\nabla\omega|^{2}dxdy+\int_{r>r_{1},0<\theta<\pi}|\nabla\boldsymbol{w}|^{2}dxdy\right).
\end{aligned}
\end{equation*}
It follows from (\ref{eq:Dirichlet-vor}) in Lemma \ref{lem:Dirichlet-vor} and (\ref{eq:Dirichlet}) that
$$ \int_{r>r_{1},0<\theta<\pi}\frac{|\nabla p|^{2}}{\log r}<\infty$$
for any $ r_{1}>\max (r_{0},1). $ By the integral theorem of the mean and Wirtinger's inequality there is an $ R_{n}\in (2^{2^{n}},2^{2^{n+1}}) $ such that
\begin{equation*}
\begin{aligned}
\log 2\int_{0}^{\pi}|p(R_{n},\theta)-\bar{p}(R_{n})|^{2}d\theta&=\int_{2^{2^{n}}}^{2^{2^{n+1}}}\frac{1}{r\log r}dr\int_{0}^{\pi}|p(r,\theta)-\bar{p}(r)|^{2}d\theta\\&\leq\int_{2^{2^{n}}}^{2^{2^{n+1}}}\int_{0}^{\pi}\frac{p_{\theta}^{2}}{r\log r}d\theta dr\\&\leq\int_{2^{2^{n}}<r<2^{2^{n+1}},0<\theta<\pi}\frac{|\nabla p|^{2}}{\log r}dxdy\to 0, \quad n\to\infty.
\end{aligned}
\end{equation*}
Thus we have (\ref{eq:pressure convergence1}).

The next lemma is the uniform square convergence of the pressure, and we used the Green function in the half-annulus.
\begin{lemma}\label{lem:pressure convergence2}
Under the assumptions of Theorem \ref{thm:pressuredecay}, we have
	\begin{equation}\label{eq:pressure convergence2}
	\begin{aligned}
	\lim_{r\to\infty}\int_{0}^{\pi}|p(r,\theta)-\bar{p}(r)|^{2}d\theta=0.
	\end{aligned}
	\end{equation}
\end{lemma}

{\bf Proof.} Using (\ref{eq:component-NS}), we find that
$$\Delta p=2(u_{x}v_{y}-u_{y}v_{x}).$$ The right member is absolutely integrable in $ \Omega_{0}. $ It also follows that
\begin{equation}\label{eq:partial-theta pressure}
\begin{aligned}
\frac{1}{r}p_{\theta}=-\omega_{r}-uv_{r}+(v+b)u_{r},
\end{aligned}
\end{equation}
and
\beno \Delta\bar{p}&=&\int_0^{\pi}\triangle p-\frac{1}{r^2}\frac{\partial^2 p}{\partial\theta^2}d\theta\\
&=&\int_0^{\pi}\triangle p d\theta-\frac{1}{r}[-\omega_{r}-uv_{r}+(v+b)u_{r}]|_0^{\pi}\\&=&\int_0^{\pi}\triangle p d\theta,
\eeno
which is also absolutely integrable due to Navier-slip boundary condition (\ref{eq:slip}) and $ bu(r,0)=bu(r,\pi) $. Hence
\begin{equation}\label{eq:H-bound}
\begin{aligned}
H\equiv\Delta(p-\bar{p})\in L_{1}, \quad in~~\Omega_{0}.
\end{aligned}
\end{equation}

Let $ A_{nm}^{+} $ denote the upper half annulus $ R_{n}<r<R_{m},0<\theta<\pi, $ the sequence of radii $ R_{n} $ being defined as in Lemma \ref{lem:pressure convergence1}. Writing $ r=R_{n}\bar{r}, $ then
$1<\bar{r}<\frac{R_{m}}{R_{n}}\triangleq R_{0}.$ Noting that
$$\Delta[p(R_{n}\bar{r},\theta)-\bar{p}(R_{n}\bar{r})]=R_{n}^{2}\Delta[p(r,\theta)-\bar{p}(r)]=R_{n}^{2}H(r,\theta).$$
 We have the representation
\begin{equation}\label{eq:Green representation}
\begin{aligned}
p(r,\theta)-\bar{p}(r)=&p(R_{n}\bar{r},\theta)-\bar{p}(R_{n}\bar{r})\\=&-\int_{1<\bar{r}<R_{0},0<\theta<\pi}G(\bar{r},\theta;\bar{\rho},\varphi)R_{n}^{2}H(R_{n}\bar{\rho},\varphi)\bar{\rho}d\bar{\rho} d\varphi\\&+\int_{\bar{\rho}=1,0<\theta<\pi}\frac{\partial G}{\partial\rho}(\bar{r},\theta;\bar{\rho},\varphi)(p(R_{n}\bar{\rho},\varphi)-\bar{p}(R_{n}\bar{\rho}))\bar{\rho} d\varphi\\&-\int_{1<\bar{r}<R_{0},\theta=0}\frac{1}{\bar{\rho}}\frac{\partial G}{\partial\varphi}(\bar{r},\theta;\bar{\rho},\varphi)(p(R_{n}\bar{\rho},\varphi)-\bar{p}(R_{n}\bar{\rho}))\bar{\rho} d\varphi\\&+\int_{1<\bar{r}<R_{0},\theta=\pi}\frac{1}{\bar{\rho}}\frac{\partial G}{\partial\varphi}(\bar{r},\theta;\bar{\rho},\varphi)(p(R_{n}\bar{\rho},\varphi)-\bar{p}(R_{n}\bar{\rho}))\bar{\rho} d\varphi
\\&-\int_{\bar{\rho}=R_{0},0<\theta<\pi}\frac{\partial G}{\partial\rho}(\bar{r},\theta;\bar{\rho},\varphi)(p(R_{n}\bar{\rho},\varphi)-\bar{p}(R_{n}\bar{\rho}))\bar{\rho} d\varphi\\=&J_{1}+J_{2}+J_{3}+J_{4}+J_{5},
\end{aligned}
\end{equation}
where $ G=G(\bar{r},\theta;\bar{\rho},\varphi) $ is the harmonic Green's function for the upper half annulus $ 1<\bar{r}<R_{0},0<\theta<\pi.~~G $ can be written in the form (see, for example, p.140, problem 4 with answer on p.418 in \cite{Wein1965})
\begin{equation*}G(\bar{r},\theta;\bar{\rho},\varphi)= \left\{
\begin{aligned}
\sum_{k=1}^{\infty}\frac{(\bar{r}^{k}-\bar{r}^{-k})}{\pi k(R_{0}^{k}-R_{0}^{-k})}[(\frac{R_{0}}{\bar{\rho}})^{k}-(\frac{\bar{\rho}}{R_{0}})^{k}]\sin k\theta \sin k\varphi,\quad \bar{r}<\bar{\rho},\\\sum_{k=1}^{\infty}\frac{(\bar{\rho}^{k}-\bar{\rho}^{-k})}{\pi k(R_{0}^{k}-R_{0}^{-k})}[(\frac{R_{0}}{\bar{r}})^{k}-(\frac{\bar{r}}{R_{0}})^{k}]\sin k\theta \sin k\varphi,\quad \bar{r}>\bar{\rho}.
\end{aligned} \right.
\end{equation*}
Next,we use variable substitution for (\ref{eq:Green representation}). Let $ \rho=R_{n} \bar{\rho} $ and noting that $ R_{0}=\frac{R_{m}}{R_{n}},\bar{r}=\frac{r}{R_{n}}. $ For $ J_{1}, $ we have
\begin{equation*}
	\begin{aligned}
	J_{1}=&-\int_{1<\bar{r}<R_{0},0<\theta<\pi}G(\bar{r},\theta;\bar{\rho},\varphi)R_{n}^{2}H(R_{n}\bar{\rho},\varphi)\bar{\rho}d\bar{\rho} d\varphi\\=&-\int_{1<\bar{\rho}<\bar{r},0<\theta<\pi}\sum_{k=1}^{\infty}\frac{(\bar{\rho}^{k}-\bar{\rho}^{-k})}{\pi k(R_{0}^{k}-R_{0}^{-k})}[(\frac{R_{0}}{\bar{r}})^{k}-(\frac{\bar{r}}{R_{0}})^{k}]\sin k\theta \sin k\varphi R_{n}^{2}H(R_{n}\bar{\rho},\varphi)\bar{\rho}d\bar{\rho}\\&-\int_{\bar{r}<\bar{\rho}<R_{0},0<\theta<\pi}\sum_{k=1}^{\infty}\frac{(\bar{r}^{k}-\bar{r}^{-k})}{\pi k(R_{0}^{k}-R_{0}^{-k})}[(\frac{R_{0}}{\bar{\rho}})^{k}-(\frac{\bar{\rho}}{R_{0}})^{k}]\sin k\theta \sin k\varphi R_{n}^{2}H(R_{n}\bar{\rho},\varphi)\bar{\rho}d\bar{\rho}
\\=&-\int_{R_{n}<\rho<r,0<\theta<\pi}[-\sum_{k=1}^{\infty}\frac{(\rho^{k}-R_{n}^{2k}/\rho^{k})}{\pi k(R_{m}^{2k}-R_{n}^{2k})}(r^{k}-R_{m}^{2k}/r^{k})\sin k\theta\sin k\varphi]H(\rho,\varphi)\rho d\rho d\varphi\\&-\int_{r<\rho<R_{m},0<\theta<\pi}[-\sum_{k=1}^{\infty}\frac{(r^{k}-R_{n}^{2k}/r^{k})}{\pi k(R_{m}^{2k}-R_{n}^{2k})}(\rho^{k}-R_{m}^{2k}/\rho^{k})\sin k\theta\sin k\varphi]H(\rho,\varphi)\rho d\rho d\varphi.
	\end{aligned}
\end{equation*}
For $ J_{2} ,$ noting that $ \bar{\rho}=1<\bar{r}, $ we get
\begin{equation*}
	\begin{aligned}
	J_{2}&=\int_{\bar{\rho}=1,0<\theta<\pi}\frac{\partial G}{\partial\rho}(\bar{r},\theta;\bar{\rho},\varphi)(p(R_{n}\bar{\rho},\varphi)-\bar{p}(R_{n}\bar{\rho}))\bar{\rho} d\varphi\\&=\int_{\bar{\rho}=1,0<\theta<\pi}\frac{\partial\left\{\sum_{k=1}^{\infty}\frac{\bar{\rho}^{k}-\bar{\rho}^{-k}}{\pi k(R_{0}^{k}-R_{0}^{-k})}[(\frac{R_{0}}{\bar{r}})^{k}-(\frac{\bar{r}}{R_{0}})^{k}]\sin k\theta\sin k\varphi\right\}}{\partial\bar{\rho}}[p(R_{n}\bar{\rho},\varphi)-\bar{p}(R_{n}\bar{\rho})]\bar{\rho}d\varphi
\\&=\int_{\rho=R_{n},0<\theta<\pi}\frac{\partial[-\sum_{k=1}^{\infty}\frac{(\rho^{k}-R_{n}^{2k}/\rho^{k})}{\pi k(R_{m}^{2k}-R_{n}^{2k})}(r^{k}-R_{m}^{2k}/r^{k})\sin k\theta\sin k\varphi]}{\partial\rho}[p(\rho,\varphi)-\bar{p}(\rho)]\rho d\varphi.
	\end{aligned}
\end{equation*}
For $ J_{5}, $ noting that $ \bar{\rho}=R_{0}>\bar{r}, $ we obtain
\begin{equation*}
	\begin{aligned}
	J_{5}&=-\int_{\bar{\rho}=R_{0},0<\theta<\pi}\frac{\partial G}{\partial\rho}(\bar{r},\theta;\bar{\rho},\varphi)(p(R_{n}\bar{\rho},\varphi)-\bar{p}(R_{n}\bar{\rho}))\bar{\rho} d\varphi\\&=-\int_{\bar{\rho}=R_{0},0<\theta<\pi}\frac{\partial\left\{\sum_{k=1}^{\infty}\frac{\bar{r}^{k}-\bar{r}^{-k}}{\pi k(R_{0}^{k}-R_{0}^{-k})}[(\frac{R_{0}}{\bar{\rho}})^{k}-(\frac{\bar{\rho}}{R_{0}})^{k}]\sin k\theta\sin k\varphi\right\}}{\partial\bar{\rho}}[p(R_{n}\bar{\rho},\varphi)-\bar{p}(R_{n}\bar{\rho})]\bar{\rho}d\varphi\\&=-\int_{\rho=R_{m},0<\theta<\pi}\frac{\partial[-\sum_{k=1}^{\infty}\frac{(r^{k}-R_{n}^{2k}/r^{k})}{\pi k(R_{m}^{2k}-R_{n}^{2k})}(\rho^{k}-R_{m}^{2k}/\rho^{k})\sin k\theta\sin k\varphi]}{\partial\rho}[p(\rho,\varphi)-\bar{p}(\rho)]\rho d\varphi.
	\end{aligned}
\end{equation*}
Since $ G(\bar{r},\theta;\bar{\rho},\varphi)=0 $ for $\theta=0,\pi.$ It follows that $ J_{3}=J_{4}=0. $ Hence the representation of $G(r,\theta;\rho,\varphi)$ is follows:
\begin{equation*}G(r,\theta;\rho,\varphi)= \left\{
\begin{aligned}
-\sum_{k=1}^{\infty}\frac{(r^{k}-R_{n}^{2k}/r^{k})}{\pi k(R_{m}^{2k}-R_{n}^{2k})}(\rho^{k}-R_{m}^{2k}/\rho^{k})\sin k\theta\sin k\varphi,\quad r<\rho,\\-\sum_{k=1}^{\infty}\frac{(\rho^{k}-R_{n}^{2k}/\rho^{k})}{\pi k(R_{m}^{2k}-R_{n}^{2k})}(r^{k}-R_{m}^{2k}/r^{k})\sin k\theta\sin k\varphi,\quad r>\rho.
\end{aligned} \right.
\end{equation*}
Moreover, (\ref{eq:Green representation}) can be written in the following form
\begin{equation}\label{eq:Green function}
\begin{aligned}
p(r,\theta)-\bar{p}(r)&=-\int_{A_{nm}^{+}}G(r,\theta;\rho,\varphi)H(\rho,\varphi)\rho d\rho d\varphi\\&+\int_{\rho=R_{n},0<\theta<\pi}\frac{\partial G}{\partial\rho}(r,\theta;\rho,\varphi)(p(\rho,\varphi)-\bar{p}(\rho))\rho d\varphi\\&-\int_{\rho=R_{m},0<\theta<\pi}\frac{\partial G}{\partial\rho}(r,\theta;\rho,\varphi)(p(\rho,\varphi)-\bar{p}(\rho))\rho d\varphi.
\end{aligned}
\end{equation}
Define
\begin{equation*}
\begin{aligned}
&G^{(2)}(r;\rho_{1},\varphi_{1};\rho_{2},\varphi_{2})\\&=\int_{0}^{\pi}G(r,\theta;\rho_{1},\varphi_{1})G(r,\theta;\rho_{2},\varphi_{2})d\theta
\\&=\sum_{k=1}^{\infty}\frac{(r^{k}-R_{n}^{2k}/r^{k})^{2}(\rho_{1}^{k}-R_{m}^{2k}/\rho_{1}^{k})(\rho_{2}^{k}-R_{m}^{2k}/\rho_{2}^{k})\sin k\varphi_{1}\sin k\varphi_{2}}{2\pi k^{2}(R_{m}^{2k}-R_{n}^{2k})^{2}}
\end{aligned}
\end{equation*}
when $ r<\rho_{1},r<\rho_{2}, $ with similar expressions for the other cases. Therefore
\begin{equation*}
\begin{aligned}
|G^{(2)}(r;\rho_{1},\varphi_{1};\rho_{2},\varphi_{2})|&\leq\sum_{k=1}^{\infty}\frac{(r^{k}-R_{n}^{2k}/r^{k})^{2}|\rho_{1}^{k}-R_{m}^{2k}/\rho_{1}^{k}||\rho_{2}^{k}-R_{m}^{2k}/\rho_{2}^{k}|}{2\pi k^{2}(R_{m}^{2k}-R_{n}^{2k})^{2}}\\&:=\hat{G}^{(2)}.
\end{aligned}
\end{equation*}
Due to $ r>R_{n}, $ we find $ (r^{k}-R_{n}^{2k}/r^{k})^{2} $ is monotonically increasing with respect to $ r. $ When $ r=\rho_{1} $ or $ r=\rho_{2} ,$ $ \hat{G}^{(2)} $ attains at the maximum. If $ \rho_{1}<\rho_{2}, $ then
$$\hat{G}^{(2)}\leq \sum_{k=1}^{\infty}\frac{(\rho_{1}^{k}-R_{n}^{2k}/\rho_{1}^{k})^{2}|\rho_{1}^{k}-R_{m}^{2k}/\rho_{1}^{k}||\rho_{2}^{k}-R_{m}^{2k}/\rho_{2}^{k}|}{2\pi k^{2}(R_{m}^{2k}-R_{n}^{2k})^{2}}.$$
Noting that $ |\rho_{2}^{k}-R_{m}^{2k}/\rho_{2}^{k}| $ is monotonically decreasing with respect to $ \rho_{2}, $  $ \hat{G}^{(2)} $ reaches its maximum at $ \rho_{1}=\rho_{2}. $ Let
$$ h(r)=(r^{k}-\frac{R_{n}^{2k}}{r^{k}})(r^{k}-\frac{R_{m}^{2k}}{r^{k}}).$$ An easy computation which lets $ h'(r)=0 $ shows that $ r=(R_{m}R_{n})^{\frac{1}{2}}. $
In conclusion, the maximum with respect to $\rho_{1}$ and $\rho_{2}$ of $ \hat{G}^{(2)} $ occurs when $\rho_{1}=\rho_{2}=r=(R_{m}R_{n})^{\frac{1}{2}}.$ Thus
\begin{equation*}
\begin{aligned}
|G^{(2)}(r;\rho_{1},\varphi_{1};\rho_{2},\varphi_{2})|&\leq\sum_{k=1}^{\infty}\frac{[(R_{m}R_{n})^{\frac{k}{2}}-R_{n}^{2k}/(R_{m}R_{n})^{\frac{k}{2}}]^{2}[(R_{m}R_{n})^{\frac{k}{2}}-R_{m}^{2k}/(R_{m}R_{n})^{\frac{k}{2}}]^{2}}{2\pi k^{2}(R_{m}^{2k}-R_{n}^{2k})^{2}}\\&=\sum_{k=1}^{\infty}\frac{[R_{m}^{-\frac{k}{2}}R_{n}^{\frac{k}{2}}(R_{m}^{k}-R_{n}^{k})]^{2}[R_{m}^{\frac{k}{2}}R_{n}^{-\frac{k}{2}}(R_{n}^{k}-R_{m}^{k})]^{2}}{2\pi k^{2}(R_{m}^{k}+R_{n}^{k})^{2}(R_{m}^{k}-R_{n}^{k})^{2}}
\\&=\sum_{k=1}^{\infty}\frac{(R_{m}^{k}-R_{n}^{k})^{2}}{2\pi k^{2}(R_{m}^{k}+R_{n}^{k})^{2}}\leq\sum_{k=1}^{\infty}\frac{1}{k^{2}}\equiv C_{1}.
\end{aligned}
\end{equation*}
Also
\begin{equation*}
\begin{aligned}
\frac{\partial G(r,\theta;\rho,\varphi)}{\partial\rho}|_{\rho=R_{m}}
&=\frac{\partial[-\sum_{k=1}^{\infty}\frac{(r^{k}-R_{n}^{2k}/r^{k})}{\pi k(R_{m}^{2k}-R_{n}^{2k})}(\rho^{k}-R_{m}^{2k}/\rho^{k})\sin k\theta\sin k\varphi]}{\partial\rho}|_{\rho=R_{m}}\\&=-\sum_{k=1}^{\infty}\frac{(r^{k}-R_{n}^{2k}/r^{k})(k\rho^{k-1}+kR_{m}^{2k}/\rho^{k+1})}{\pi k(R_{m}^{2k}-R_{n}^{2k})}\sin k\theta\sin k\varphi|_{\rho=R_{m}}\\&=
-\sum_{k=1}^{\infty}\frac{2(r^{k}-R_{n}^{2k}/r^{k})R_{m}^{k-1}}{\pi(R_{m}^{2k}-R_{n}^{2k})}\sin k\theta\sin k\varphi
\end{aligned}
\end{equation*}
and thus, if $ R_{n}<r\leq R_{m-2},0<\theta<\pi, $ we have
\begin{equation*}
\begin{aligned}
&\int_{0}^{\pi}|G_{\rho}(r,\theta;R_{m},\varphi)|^{2}R_{m}^{2}d\varphi
\\&=\sum_{k=1}^{\infty}\frac{4(r^{k}-R_{n}^{2k}/r^{k})^{2}R_{m}^{2k}}{\pi^{2}(R_{m}^{2k}-R_{n}^{2k})^{2}}\sin^{2}k\theta\int_{0}^{\pi}\sin^{2}k\varphi d\varphi\\&\leq\sum_{k=1}^{\infty}\frac{2(r^{k}-R_{n}^{2k}/r^{k})^{2}R_{m}^{2k}}{\pi(R_{m}^{2k}-R_{n}^{2k})^{2}}
\\&\leq \sum_{k=1}^{\infty}\frac{2(R_{m-2}/R_{m})^{2k}}{\pi[1-(R_{n}/R_{m})^{2k}]^{2}}\\&\leq\frac{\sum_{k=1}^{\infty}2^{-2^{m}k}}{[1-(R_{n}/R_{m})^{2}]^{2}}\leq C_{2}.
\end{aligned}
\end{equation*}
Similarly,  if $ R_{n+2}\leq r<R_{m},0<\theta<\pi $, we have
$$\frac{\partial G(r,\theta;\rho,\varphi)}{\partial\rho}|_{\rho=R_{n}}=-\sum_{k=1}^{\infty}\frac{2(r^{k}-R_{m}^{2k}/r^{k})R_{n}^{k-1}}{\pi(R_{m}^{2k}-R_{n}^{2k})}\sin k\theta\sin k\varphi$$
and
\begin{equation*}
	\begin{aligned}
	&\int_{0}^{\pi}|G_{\rho}(r,\theta;R_{n},\varphi)|^{2}R_{n}^{2}d\varphi\\&=\sum_{k=1}^{\infty}\frac{4(r^{k}-R_{m}^{2k}/r^{k})^{2}R_{n}^{2k}}{\pi^{2}(R_{m}^{2k}-R_{n}^{2k})^{2}}\sin^{2}k\theta\int_{0}^{\pi}\sin^{2}k\varphi d\varphi\\&\leq\sum_{k=1}^{\infty}\frac{2(R_{n}/R_{n+2})^{2k}}{\pi[1-(R_{n}/R_{m})^{2k}]^{2}}\\&\leq \frac{\sum_{k=1}^{\infty}2^{-2^{n}k}}{[1-(R_{n}/R_{m})^{2}]^{2}}\leq C_{3}.
	\end{aligned}
\end{equation*}
By H$ \ddot{\rm{o}} $lder inequality, it follows from (\ref{eq:Green function}) that for $ r\in [R_{n+2},R_{m-2}],0<\theta<\pi,m\geq n+5 $
\begin{equation*}
\begin{aligned}
&\int_{0}^{\pi}|p(r,\theta)-\bar{p}(r)|^{2}d\theta\\\leq& 3\int_{0}^{\pi}[\int_{0}^{\pi}\int_{R_{n}}^{R_{m}}G(r,\theta;\rho,\varphi)H(\rho,\varphi)\rho d\rho d\varphi]^{2}d\theta\\&+3\int_{0}^{\pi}[\int_{0}^{\pi}\frac{\partial G}{\partial\rho}(r,\theta;R_{n},\varphi)(p(R_{n},\varphi)-\bar{p}(R_{n}))R_{n}d\varphi]^{2}d\theta\\&+3\int_{0}^{\pi}[\int_{0}^{\pi}\frac{\partial G}{\partial\rho}(r,\theta;R_{m},\varphi)(p(R_{m},\varphi)-\bar{p}(R_{m}))R_{m}d\varphi]^{2}d\theta\\\leq& 3\int_{0}^{\pi}\int_{R_{n}}^{R_{m}}\int_{0}^{\pi}\int_{R_{n}}^{R_{m}}[\int_{0}^{\pi}G(r,\theta;\rho_{1},\varphi_{1})G(r,\theta;\rho_{2},\varphi_{2})d\theta] H(\rho_{1},\varphi_{1})H(\rho_{2},\varphi_{2})\rho_{1}\rho_{2}d\rho_{1}d\varphi_{1} d\rho_{2}d\varphi_{2}\\&+3\int_{0}^{\pi}[\int_{0}^{\pi}|\frac{\partial G}{\partial\rho}(r,\theta;R_{n},\varphi)|^{2}R_{n}^{2}d\varphi\cdot\int_{0}^{\pi}|p(R_{n},\varphi)-\bar{p}(R_{n})|^{2}d\varphi]d\theta\\&+3\int_{0}^{\pi}[\int_{0}^{\pi}|\frac{\partial G}{\partial\rho}(r,\theta;R_{m},\varphi)|^{2}R_{m}^{2}d\varphi\cdot\int_{0}^{\pi}|p(R_{m},\varphi)-\bar{p}(R_{m})|^{2}d\varphi]d\theta\\\leq& 3\pi C_{1}(\int_{R_{n}<r<R_{m},0<\theta<\pi}|H|dxdy)^{2}+3\pi C_{3}\int_{0}^{\pi}|p(R_{n},\theta)-\bar{p}(R_{n})|^{2}d\varphi\\&+3\pi C_{2}\int_{0}^{\pi}|p(R_{m},\varphi)-\bar{p}(R_{m})|^{2}d\varphi.
\end{aligned}
\end{equation*}
By letting $ m\to\infty $ and using (\ref{eq:pressure convergence1}) and (\ref{eq:H-bound}), we obtain an upper bound on the left member for $ r>2^{2^{n+2}}, $ and this bound approaches zero as $ n\to\infty. $ From this we infer
$$\lim_{r\to\infty}\int_{0}^{\pi}|p(r,\theta)-\bar{p}(r)|^{2}d\theta=0.$$

The third lemma is about the convergence of the average pressure.
\begin{lemma}\label{lem:meanpressure convergence}
	Under the assumptions of Theorem \ref{thm:pressuredecay}, the average pressure
	$$\bar{p}(r)=\frac{1}{\pi}\int_{0}^{\pi}p(r,\theta)d\theta$$
	has a limit at infinity
	\begin{equation}\label{eq:meanpressure convergence}
	\lim_{r\to\infty}\bar{p}(r)=p_{\infty}<\infty.
	\end{equation}
\end{lemma}

{\bf Proof.} First,
 Navier-slip boundary condition tells us that $ \omega(x,0)=0 $ and $ v(x,0)=0 $ due to (\ref{eq:slip}). Hence
 \begin{equation*}
 	\int_{0}^{\pi}\omega_{\theta}d\theta=0,
 \end{equation*}
 and
 \begin{equation*}
 	\int_{0}^{\pi}(uv_{\theta}+vu_{\theta})d\theta=\int_{0}^{\pi}\bar{u}v_{\theta}d\theta=0,
 \end{equation*}
where
$$ \bar{u}(r)=\frac{1}{\pi}\int_{0}^{\pi}u(r,\theta)d\theta.$$
Besides, due to $ bu(r,0)=bu(r,\pi) $,  we average (\ref{eq:partial pressure}) to find that
\begin{equation}\label{eq:mean partial pressure}
	\begin{aligned}
	\bar{p}'(r)&=\frac{1}{\pi r}\int_{0}^{\pi}[\omega_{\theta}+uv_{\theta}-(v+b)u_{\theta}]d\theta\\&=\frac{1}{\pi r}\int_{0}^{\pi}[2uv_{\theta}-(vu_{\theta}+uv_{\theta})]d\theta
\\&=\frac{2}{\pi r}\int_{0}^{\pi}(u-\bar{u})v_{\theta}d\theta.
	\end{aligned}
\end{equation}
Integrating this inequality with respect to $ r $ over $ (r_{1},r_{2})~~ (r_{2}\geq r_{1}\geq r_{0}), $ then by Cauchy and Wirtinger inequalities, we find
\begin{equation*}
\begin{aligned}
|\bar{p}(r_{2})-\bar{p}(r_{1})|&=|\frac{2}{\pi}\int_{r_{1}}^{r_{2}}\int_{0}^{\pi}\frac{(u-\bar{u})v_{\theta}}{r}d\theta dr|
\\&\leq \frac{1}{\pi}\int_{r_{1}}^{r_{2}}\int_{0}^{\pi}\frac{|\boldsymbol{w}-\bar{\boldsymbol{w}}|^{2}+|\boldsymbol{w}_{\theta}|^{2}}{r}d\theta dr
\\&\leq C\int_{r_{1}}^{r_{2}}\int_{0}^{\pi}\frac{|\boldsymbol{w}_{\theta}|^{2}}{r}d\theta dr\\&\leq C\int_{r>r_{1}}|\nabla\boldsymbol{w}|^{2}dxdy.
\end{aligned}
\end{equation*}
Since the right member of this inequality tends to zero as $ r_{1}\to\infty, $ it follows that $ \bar{p}(r) $ has a limit $ p_{\infty}, $ as asserted.

Immediately it follows from Lemma \ref{lem:meanpressure convergence} and Lemma \ref{lem:pressure convergence2} that the following conclusion holds.
\begin{corollary}\label{coro:pressure convergence}
	Under the assumptions of Theorem \ref{thm:pressuredecay}, we have
	\ben\label{eq:pressure convergence3}
	\lim_{r\to\infty}\int_{0}^{\pi}|p(r,\theta)-p_{\infty}|^{2}d\theta=0.
	\een
\end{corollary}


{\bf Proof of Case (i) in Theorem \ref{thm:pressuredecay}.} For $ r>r_{0}, $ define
$$\tilde{p}(\tilde{x})=r^{2}p(r\tilde{x})=r^{2}p(x),$$
where $ x\in B_{2r}^{+}\backslash B_{r}^{+} $ and $ r $ is large enough. Without loss of generality, we still consider $B_{r}^{+}$, since one can mollify the domain such that it's regular. By Lemma \ref{lem:BGI}, we have
\begin{equation*}
\begin{aligned}
\left\|\tilde{p}\right\|_{L^{\infty}(B_{\frac74}^{+}\backslash B_{\frac54}^{+})}&\leq C(1+\left\|\tilde{p}\right\|_{H^{1}(B_{2}^{+}\backslash B_{1}^{+})})\sqrt{\log (e+\left\|\Delta\tilde{p}\right\|_{L^{2}(B_{2}^{+}\backslash B _{1}^{+})})}\\&=C(1+\left\|\tilde{p}\right\|_{L^{2}(B_{2}^{+}\backslash B_{1}^{+})}+\left\|\nabla \tilde{p}\right\|_{L^{2}(B_{2}^{+}\backslash B_{1}^{+})})\sqrt{\log (e+\left\|\Delta\tilde{p}\right\|_{L^{2}(B_{2}^{+}\backslash B _{1}^{+})})}.
\end{aligned}
\end{equation*}
Due to the scaling, we have
$$\left\|\tilde{p}\right\|_{L^{\infty}(B_{2}^{+}\backslash B_{1}^{+})}=r^{2}\left\|p\right\|_{L^{\infty}(B_{2r}^{+}\backslash B_{r}^{+})},$$
$$\left\|\tilde{p}\right\|_{L^{2}(B_{2}^{+}\backslash B_{1}^{+})}=r\left\|p\right\|_{L^{2}(B_{2r}^{+}\backslash B_{r}^{+})},$$
$$\left\|\nabla\tilde{p}\right\|_{L^{2}(B_{2}^{+}\backslash B_{1}^{+})}=r^{2}\left\|\nabla p\right\|_{L^{2}(B_{2r}^{+}\backslash B_{r}^{+})},$$
$$\left\|\nabla^{2}\tilde{p}\right\|_{L^{2}(B_{2}^{+}\backslash B_{1}^{+})}=r^{3}\left\|\nabla^{2} p\right\|_{L^{2}(B_{2r}^{+}\backslash B_{r}^{+})}.$$
Hence
\begin{equation}\label{eq:BGI}
\begin{aligned}
r^{2}\left\|p\right\|_{L^{\infty}(B_{\frac74 r}^{+}\backslash B_{\frac54 r}^{+})}\leq&C(1+r\left\|p\right\|_{L^{2}(B_{2r}^{+}\backslash B_{r}^{+})}+r^{2}\left\|\nabla p\right\|_{L^{2}(B_{2r}^{+}\backslash B_{r}^{+})})\\&\cdot\sqrt{\log (e+r^{3}\left\|\Delta p\right\|_{L^{2}(B_{2r}^{+}\backslash B_{r}^{+})})}.
\end{aligned}
\end{equation}
From Lemma \ref{lem:meanpressure convergence} we know $ |p_{\infty}|\leq C $ and by (\ref{eq:pressure convergence3}) in Corollary \ref{coro:pressure convergence} we have
\begin{equation*}
\begin{aligned}
\left\|p\right\|_{L^{2}(B_{2r}^{+}\backslash B_{r}^{+})}^{2}&\leq 2\int_{r}^{2r}\int_{0}^{\pi}|p(\rho,\theta)-p_{\infty}|^{2}d\theta\rho d\rho+2\int_{r}^{2r}\int_{0}^{\pi}|p_{\infty}|^{2}\rho d\theta d\rho\\&\leq  o(r)\int_{r}^{2r}\rho d\rho+Cr^{2}.
\end{aligned}
\end{equation*}
Then
\begin{equation*}
\left\|p\right\|_{L^{2}(B_{2r}^{+}\backslash B_{r}^{+})}^{2}\leq Cr^{2}.
\end{equation*}
Since $ |\Delta\boldsymbol{w}|^{2}=|\nabla\omega|^{2}, $ then by Lemma \ref{lem:Dirichlet-vor}, $ \boldsymbol{w}=o(\sqrt{\log r}) $ and Navier-Stokes equation $ (\ref{eq:NS})_{1}, $ there holds
\begin{equation*}
\begin{aligned}
&\left\|\nabla p\right\|_{L^{2}(B_{2r}^{+}\backslash B_{r}^{+})}^{2}\\&\leq C\int_{B_{2r}^{+}\backslash B_{r}^{+}}|\Delta\boldsymbol{w}|^{2}dxdy+C\int_{B_{2r}^{+}\backslash B_{r}^{+}}|\boldsymbol{w}\cdot\nabla\boldsymbol{w}|^{2}dxdy+C\int_{B_{2r}^{+}\backslash B_{r}^{+}}|\nabla\boldsymbol{w}|^{2}dxdy\\&\leq C\int_{B_{2r}^{+}\backslash B_{r}^{+}}|\nabla\omega|^{2}dxdy+C[o(\log r)+1]\int_{B_{2r}^{+}\backslash B_{r}^{+}}|\nabla\boldsymbol{w}|^{2}dxdy\\&\leq C[1+o(\log r)].	
\end{aligned}
\end{equation*}
Besides, similar to the second step of Theorem \ref{thm:velocitydecay}, using (\ref{eq:extension-vel1}) and (\ref{eq:extension-vel2}),we also get
\begin{equation*}
\begin{aligned}
\|\nabla^{2}(\hat{\boldsymbol{w}_{\varepsilon}}\phi)\|_{L^{2}(\R^{2})}&\leq C\|\nabla(\hat{\omega_\varepsilon}\phi+\hat{\boldsymbol{w}_\varepsilon}\cdot\nabla^{\bot}\phi)\|_{L^{2}(\R^{2})}+C\|\nabla[\nabla\cdot(\hat{\boldsymbol{w}_{\varepsilon}}\phi)]\|_{L^{2}(\mathbb{R}^{2})}\\&=C\|\phi\nabla\hat{\omega_\varepsilon}+\hat{\omega_\varepsilon}\nabla\phi+\nabla(\hat{\boldsymbol{w}_{\varepsilon}}\cdot\nabla^{\bot}\phi)\|_{L^{2}(\mathbb{R}^{2})}+C\|\nabla[(\nabla\phi)\hat{\boldsymbol{w}_{\varepsilon}}]\|_{L^{2}(\mathbb{R}^{2})}\\&\leq C\|\nabla\hat{\omega_\varepsilon}\|_{L^{2}(\R^{2}\backslash B_{2r_{0}})}+C\|\hat{\omega_\varepsilon}\|_{L^{2}(B_{3r_{0}}\backslash B_{2r_{0}})}+C\|\nabla\hat{\boldsymbol{w}_\varepsilon}\|_{L^{2}(B_{3r_{0}\backslash }B_{2r_{0}})}.
\end{aligned}
\end{equation*}
Since $ \hat{\omega_\varepsilon}\in C^{\infty}(\R^{2}\backslash B_{2r_{0}}), $ we know $ \|\hat{\omega_\varepsilon}\|_{L^{2}(B_{3r_{0}}\backslash B_{2r_{0}})}\leq C. $ Using (\ref{eq:Dirichlet}), we have $ \|\nabla\hat{\boldsymbol{w}_\varepsilon}\|_{L^{2}(B_{3r_{0}\backslash }B_{2r_{0}})}\leq C. $ Therefore
\begin{equation*}
\left\|\nabla^{2}\hat{\boldsymbol{w}_{\varepsilon}}\right\|_{L^{2}(\R^{2}\backslash B_{3r_{0}})}\leq C+C\|\nabla\hat{\omega_\varepsilon}\|_{L^{2}(\R^{2}\backslash B_{2r_{0}})}\leq C+C\|\nabla\hat{\omega}\|_{L^{2}(\R^{2}\backslash B_{2r_{0}})}.
\end{equation*}
Let $\varepsilon\to 0,$ by Lebesgue's dominated convergence theorem and Lemma \ref{lem:Dirichlet-vor}, we have
\begin{equation*}
\begin{aligned}
\|\nabla^{2}\hat{\boldsymbol{w}}\|_{L^{2}(\R^{2}_+\backslash B_{3r_{0}})}&\leq C+C\|\nabla\hat{\omega}\|_{L^{2}(\R^{2}\backslash B_{2r_{0}})}\\&\leq C+C\|\nabla\omega\|_{L^{2}(\R_{+}^{2}\backslash B_{2r_{0}}^{+})}\leq C.
\end{aligned}
\end{equation*}
Then by $ (\ref{eq:component-NS})_{1,2} $ and Gagliardo-Nirenberg inequality, we have
\begin{equation*}
\begin{aligned}
&\int_{B_{2r}^{+}\backslash B_{r}^{+}}|\Delta p|^{2}dxdy\\\leq &C\int_{B_{2r}^{+}\backslash B_{r}^{+}}(u_{x}^{2}v_{y}^{2}+u_{y}^{2}v_{x}^{2})dxdy\\\leq& C\int_{B_{2r}^{+}\backslash B_{r}^{+}}|\nabla\boldsymbol{w}|^{4}dxdy\\\leq& C(\int_{B_{2r}^{+}\backslash B_{r}^{+}}|\nabla\boldsymbol{w}|^{2}dxdy)(\int_{B_{2r}^{+}\backslash B_{r}^{+}}|\nabla^{2}\boldsymbol{w}|^{2} dxdy)+Cr^{-2}(\int_{B_{2r}^{+}\backslash B_{r}^{+}}|\nabla\boldsymbol{w}|^{2}dxdy)^{2}\\\leq &C(1+r^{-2}).
\end{aligned}
\end{equation*}
Therefore using (\ref{eq:BGI}), one can get
\begin{equation*}
\begin{aligned}
r^{2}\left\|p\right\|_{L^{\infty}(B_{\frac74r}^{+}\backslash B_{\frac54r}^{+})}&\leq C(1+Cr^{2}+Cr^{2}\sqrt{1+o(\log r)})\sqrt{\log(e+Cr^{3}\sqrt{1+r^{-2}})}\\&\leq r^{2}o(\log r),
\end{aligned}
\end{equation*}
which implies
$$\left\|p\right\|_{L^{\infty}(B_{\frac74r}^{+}\backslash B_{\frac54r}^{+})}\leq o(\log r).$$
The proof is complete.

{\bf  Proof of Case (ii) in Theorem \ref{thm:pressuredecay}.} Let the point $ P(2R,\theta) $ be the origin of a new system of polar coordinates $ (r',\theta') $ and suppose that $ R>r_{0}. $ In these new coordinates we still have
\begin{equation*}
	p_{r'}=\frac{1}{r'}[\omega_{\theta'}+uv_{\theta'}-(v+b)u_{\theta'}]
\end{equation*}
from the Navier-Stokes equations.
Integrating with respect to $ r' $ over $ (0,r') $
\begin{equation}\label{eq:pressure expression}
	p(P)=p(r',\theta')+\int_{0}^{r'}\frac{1}{\rho}[(v+b)u_{\theta'}-uv_{\theta'}-\omega_{\theta'}]d\rho.
\end{equation}

$\mathbf{Proof~of~Case~\uppercase\expandafter{\romannumeral1}:}$ $ p(2R,\theta) $ with $\theta=0, \pi.$ As shown in Figure 1:
\begin{figure}[H]
		\centering
\begin{tikzpicture}
	\draw[->,line width=1.0pt](-5,0)--(5,0) node[right]{$x$};
	\draw[->,line width=1.0pt](0,0)--(0,5) node[right]{$y$};
	\filldraw[black](0,0) circle (2pt);
	\draw(0,0)--(0,0.1) node[below=2pt]{$0$};
	\draw(1,0)--(1,0.1) node[below=2pt]{$R$};
	\draw(2,0)--(2,0.1) node[below=2pt]{$2R$};
	\draw(3,0)--(3,0.1) node[below=2pt]{$3R$};
	\draw(4,0)--(4,0.1) node[below=2pt]{$4R$};
	\draw(-1,0)--(-1,0.1) node[below=2pt]{$-R$};
	\draw(-2,0)--(-2,0.1) node[below=2pt]{$-2R$};
	\draw(-3,0)--(-3,0.1) node[below=2pt]{$-3R$};
	\draw(-4,0)--(-4,0.1) node[below=2pt]{$-4R$};
	\draw(0,1)--(0.1,1) node[left=2pt]{$R$};
	\draw(0,2)--(0.1,2) node[left=2pt]{$2R$};
	\draw(0,3)--(0.1,3) node[left=2pt]{$3R$};
	\draw(0,4)--(0.1,4) node[left=2pt]{$4R$};
	\draw [red,line width=1.0pt] (3,0) arc(0:180:3);
	\draw [red,line width=1.0pt] (1,0) arc(0:180:1);
	\draw [blue,line width=1.0pt] (3,0) arc(0:180:1);
	\draw [blue,line width=1.0pt] (-1,0) arc(0:180:1);
	\filldraw[blue] (2,0) circle (2pt) node [above=1pt]{$ P_{1}(2R,0) $};
	\filldraw[blue] (-2,0) circle (2pt) node [above=1pt]{$ P_{2}(2R,\pi) $};
\end{tikzpicture}
\caption{The case of $ \theta=0,\pi. $}
\end{figure}
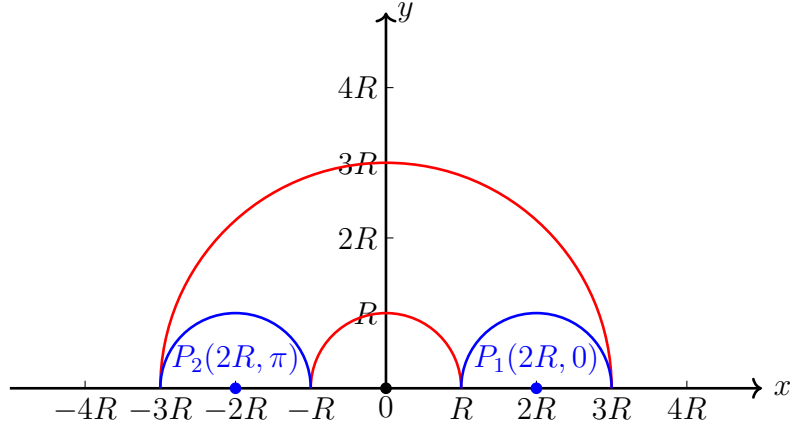
 let's consider the case of $\theta=0,$ and $\theta=\pi$ is similar.
Integrating (\ref{eq:pressure expression}) with respect to $ \theta' $ over $ [0,\pi] $
$$\pi p(P)=\int_{0}^{\pi}p(r',\theta')d\theta'+\int_{0}^{r'}\int_{0}^{\pi}\frac{1}{\rho}[(v+b)u_{\theta'}-uv_{\theta'}-\omega_{\theta'}]d\rho d\theta'.$$ Due to $ \omega(x,0)=0 $,$ v(x,0)=0 $ and $ bu(r,0)=bu(r,\pi) $. Similar to the (\ref{eq:mean partial pressure}) calculation
\begin{equation*}
\begin{aligned}
p(P)=\frac{1}{\pi}\int_{0}^{\pi}p(r',\theta')d\theta'+\frac{2}{\pi}\int_{0}^{r'}\int_{0}^{\pi}\frac{1}{\rho}[\tilde{u}(\rho)-u(\rho,\theta')]v_{\theta'}(\rho,\theta')d\rho d\theta',
\end{aligned}
\end{equation*}
where
$$ \tilde{u}(r')=\frac{1}{\pi}\int_{0}^{\pi}u(r',\theta')d\theta'.$$
Multiply this relation by $ r' $ and integrate from $ 0 $ to $ R $, and we  find
\begin{equation}\label{eq:pressure expression1}
\begin{aligned}
p(P)&=\frac{2}{\pi R^{2}}\int_{0}^{R}\int_{0}^{\pi}p(r',\theta')r'dr'd\theta'\\&+\frac{4}{\pi R^{2}}\int_{0}^{R}\int_{0}^{r'}\int_{0}^{\pi}\frac{(\tilde{u}-u)v_{\theta'}}{\rho}r'd\theta'd\rho dr' \doteq I_{1}+I_{2}.
\end{aligned}
\end{equation}
Note that the  upper half disc $ r'<R,0<\theta'<\pi $ is contained in the upper half annulus $ R<r<3R,0<\theta<\pi .$ For $ I_{1}, $ using Schwarz inequality and (\ref{eq:pressure convergence3})
\begin{equation*}
\begin{aligned}
|I_{1}|^{2}&\leq (\frac{2}{\pi R^{2}})^{2}(\int_{0}^{R}\int_{0}^{\pi}p^{2}r'dr'd\theta')(\int_{0}^{R}\int_{0}^{\pi}r'dr'd\theta')\\
&\leq
\frac{2}{\pi R^{2}}\int_{R}^{3R}rdr\left\{\max_{R<r<3R}\int_{0}^{\pi}p(r,\theta)^{2}d\theta \right\}\\&=\frac{8}{\pi}\max_{R<r<3R}\int_{0}^{\pi}p(r,\theta)^{2}d\theta\to 0,\quad {\rm as}~~R\to\infty.
\end{aligned}
\end{equation*}
For $ I_{2} ,$ using Cauchy and Wirtinger inequalities
\begin{equation*}
\begin{aligned}
|\int_{0}^{\pi}(\tilde{u}-u)v_{\theta'}d\theta'|&\leq \int_{0}^{\pi}\frac{|u-\tilde{u}|^{2}+|v_{\theta'}|^{2}}{2}d\theta'\\&\leq\int_{0}^{\pi}\frac{|\boldsymbol{w}-\tilde{\boldsymbol{w}}|^{2}
+|\boldsymbol{w}_{\theta'}|^{2}}{2}d\theta'\\&\leq\int_{0}^{\pi}|\boldsymbol{w}_{\theta'}|^{2}d\theta'\leq \int_{0}^{\pi}\rho^{2}|\nabla\boldsymbol{w}|^{2}d\theta'.
\end{aligned}
\end{equation*}
Consequently, we have
\begin{equation*}
\begin{aligned}
I_{2}&\leq\frac{4}{\pi R^{2}}\int_{0}^{R}\int_{0}^{r'}\int_{0}^{\pi}|\nabla\boldsymbol{w}(\rho,\theta')|^{2}\rho d\rho d\theta'r'dr'\\&\leq \frac{4}{\pi R^{2}}\int_{R<r<3R,0<\theta<\pi}|\nabla\boldsymbol{w}|^{2}dxdy\int_{0}^{R}\rho d\rho\\&=\frac{2}{\pi}\int_{R<r<3R,0<\theta<\pi}|\nabla\boldsymbol{w}|^{2}dxdy\to 0,\quad {\rm as}~~R\to\infty.
\end{aligned}
\end{equation*}
It follows  from (\ref{eq:pressure expression1}) that
\begin{equation*}
\begin{aligned}
|p(2R,0)|\to 0,
\end{aligned}
\end{equation*}
 as $ R\to\infty. $ Thus the proof of Case \uppercase\expandafter{\romannumeral1} is complete.

 $\mathbf{Proof~of~Case~\uppercase\expandafter{\romannumeral2}:}$ $ P(2R,\theta) $ with $\theta\in(0,\pi).$
 As shown in Figure 2:
 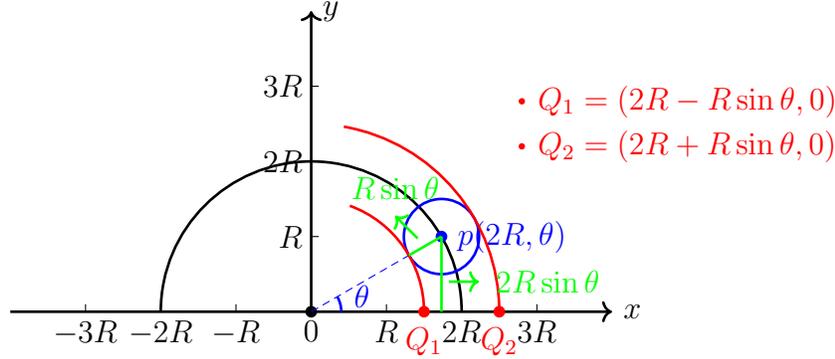
\begin{figure}[H]
 	\centering
 	\begin{tikzpicture}
 	\draw[->,line width=1.0pt](-4,0)--(4,0) node[right]{$x$};
 	\draw[->,line width=1.0pt](0,0)--(0,4) node[right]{$y$};
 	\filldraw[black](0,0) circle (2pt);
 	\draw(1,0)--(1,0.1) node[below=2pt]{$R$};
 	\draw(2,0)--(2,0.1) node[below=2pt]{$2R$};
 	\draw(3,0)--(3,0.1) node[below=2pt]{$3R$};
 	\draw(0,0)--(0,0.1) node[below=2pt]{$0$};
 	\draw(-1,0)--(-1,0.1) node[below=2pt]{$-R$};
 	\draw(-2,0)--(-2,0.1) node[below=2pt]{$-2R$};
 	\draw(-3,0)--(-3,0.1) node[below=2pt]{$-3R$};
 	\draw(0,1)--(0.1,1) node[left=2pt]{$R$};
 	\draw(0,2)--(0.1,2) node[left=2pt]{$2R$};
 	\draw(0,3)--(0.1,3) node[left=2pt]{$3R$};
 	\draw [line width=1.0pt] (2,0) arc(0:180:2);
 	\draw [blue,line width=1.0pt] (0.4,0) arc (0:30:0.4) node [right=2pt] {$ \theta $};
 	\draw[blue,densely dashed](0,0)--(1.732,1) node[right=2pt]{$ p(2R,\theta) $};
 	\filldraw[blue](1.732,1) circle (2pt);
 	\draw[blue,line width=1.0pt] (2.232,1) arc (0:360:0.5);
 	\draw[red,line width=1.0pt] (1.5,0) arc (0:70:1.5);
 	\draw[red,line width=1.0pt] (2.5,0) arc (0:80:2.5);
 	\filldraw[red] (1.5,0) circle (2pt) node [below=2pt]{$ Q_{1} $};
 	\filldraw[red] (2.5,0) circle (2pt) node [below=2pt]{$ Q_{2} $};
 	\filldraw[red] (2.8,2.8) circle (1pt) node [right=2pt] {$Q_{1}=(2R-R\sin\theta,0) $};
 	\filldraw[red] (2.8,2.2) circle (1pt) node [right=2pt] {$Q_{2}=(2R+R\sin\theta,0) $};

 	\draw[green,line width=1.0pt](1.732,1) -- (1.732,0) ;
 	\draw[->,green,line width=1.0pt](1.832,0.4)--(2.232,0.4) node [right=2pt]{$ 2R\sin\theta $};
 	\draw[green,line width=1.0pt](1.732,1)--(1.299,0.75);
 	\draw[->,green,line width=1.0pt](1.416,0.975)--(1.116,1.275) node [above=2pt] {$ R\sin\theta $};
 	\end{tikzpicture}
 	\caption{The case of $ \theta\in (0,\pi). $}
 \end{figure}

  Integrating (\ref{eq:pressure expression}) with respect to $\theta'$ over $ [0,2\pi] $
 \begin{equation*}
 	\begin{aligned}
 p(P)=\frac{1}{2\pi}\int_{0}^{2\pi}p(r',\theta')d\theta'+\frac{1}{\pi}\int_{0}^{r'}\int_{0}^{2\pi}\frac{1}{\rho}[\hat{u}(\rho)-u(\rho,\theta')]v_{\theta'}(\rho,\theta')d\rho d\theta',
\end{aligned}
\end{equation*}
where
$$ \hat{u}(r')=\frac{1}{2\pi}\int_{0}^{2\pi}u(r',\theta')d\theta'.$$
Multiply this relation by $ r' $ and integrate from $ 0 $ to $ R\sin\theta $
\begin{equation}\label{eq:pressure expression2}
\begin{aligned}
p(P)=&\frac{1}{\pi R^{2}\sin^{2}\theta}\int_{0}^{R\sin\theta}\int_{0}^{2\pi}p(r',\theta')r'dr'd\theta'\\&+\frac{2}{\pi R^{2}\sin^{2}\theta}\int_{0}^{R\sin\theta}\int_{0}^{r'}\int_{0}^{2\pi}\frac{(\hat{u}-u)v_{\theta'}}{\rho}r'd\theta'd\rho dr'\doteq& I_{1}'+I_{2}'.
\end{aligned}
\end{equation}
Noting that the disc $ r'<R\sin\theta,0<\theta'<2\pi $ is contained in the upper half annulus $ 2R-R\sin\theta<r<2R+R\sin\theta,0<\theta<\pi .$ Then similar to the Case \uppercase\expandafter{\romannumeral1}, for $ I_{1}', $ by Schwarz's inequality and (\ref{eq:pressure convergence3}), we get
\begin{equation*}
	\begin{aligned}
|I_{1}'|^{2}&\leq (\frac{1}{\pi R^{2}\sin^{2}\theta})^{2}(\int_{0}^{R\sin\theta}\int_{0}^{2\pi}p^{2}r'dr'd\theta')(\int_{0}^{R\sin\theta}\int_{0}^{2\pi}r'dr'd\theta')\\&\leq\frac{1}{\pi R^{2}\sin^{2}\theta}\int_{2R-R\sin\theta<r<2R+R\sin\theta,0<\theta<\pi}p^{2}dxdy\\&\leq
\frac{1}{\pi R^{2}\sin^{2}\theta}\int_{2R-R\sin\theta}^{2R+R\sin\theta}rdr\left\{\max_{2R-R\sin\theta<r<2R+R\sin\theta}\int_{0}^{\pi}p(r,\theta)^{2}d\theta \right\}\\&=\frac{4}{\pi\sin\theta}\max_{2R-R\sin\theta<r<2R+R\sin\theta}\int_{0}^{\pi}p(r,\theta)^{2}d\theta\to 0,\quad {\rm as}~~R\to\infty.	
	\end{aligned}
\end{equation*}
For $ I_{2}', $ using Cauchy and Wirtinger's inequalities, we have
\begin{equation*}
\begin{aligned}
I_{2}'&\leq\frac{2}{\pi R^{2}\sin^{2}\theta}\int_{0}^{R\sin\theta}\int_{0}^{r'}\int_{0}^{2\pi}|\nabla\boldsymbol{w}(\rho,\theta')|^{2}\rho d\rho d\theta'r'dr'\\&\leq \frac{2}{\pi R^{2}\sin^{2}\theta}\int_{2R-R\sin\theta<r<2R+R\sin\theta,0<\theta<\pi}|\nabla\boldsymbol{w}|^{2}dxdy\int_{0}^{R\sin\theta}\rho d\rho\\&=\frac{1}{\pi}\int_{2R-R\sin\theta<r<2R+R\sin\theta,0<\theta<\pi}|\nabla\boldsymbol{w}|^{2}dxdy\to 0,\quad {\rm as}~~R\to\infty.
\end{aligned}
\end{equation*}
Hence from (\ref{eq:pressure expression2}), we obtain $ \forall \theta\in(0,\pi), $ there holds
$$|P(2R,\theta)|\to 0$$ as $ R\to\infty. $ Hence we prove the case of \uppercase\expandafter{\romannumeral2}. Combining the results of Case \uppercase\expandafter{\romannumeral1} and Case \uppercase\expandafter{\romannumeral2} , the proof is complete.



\section{Decay  of the vorticity.}

\begin{lemma}\label{lem:vor estimate}
	Under the assumptions of Theorem \ref{thm:vorticitydecay}, we have
	\begin{equation}\label{eq:vor estimate}
	\int_{r>r_{1},0<\theta<\pi}\frac{r}{(\log r)^{\frac{1}{2}}}|\nabla\omega|^{2}dxdy<\infty(r_{1}>\max(r_{0},2-r_{0})).
	\end{equation}
\end{lemma}

{\bf Proof.} Choose $ R>r_{1}>\max(r_{0},2-r_{0}) $ and two non-negative $ C^{2} $ cut-off functions $ \xi_{1} $ and $ \xi_{2} $ such that
\begin{equation}\label{eq:cut-off function1}
\xi_{1}(r)=\left\{
\begin{aligned}
&0,\quad r\leq\frac{1}{2}(r_{0}+r_{1})\\ &1,\quad r\geq r_{1}
\end{aligned}
\right.,
\quad \xi_{2}(r)=\left\{\begin{aligned}
1,\quad r\leq 1\\ 0,\quad r\geq 2
\end{aligned}
\right..
\end{equation}
Let
$$ \eta(r)=\xi_{1}(r)\xi_{2}(\frac{r}{R})\frac{r}{(\log r)^{\frac{1}{2}}},$$
$$ h(\omega)=\omega^{2}.$$
Clearly $ \eta(r) $ vanishes near $ r=r_{0} $ and near $ r=\infty. $ Noting that $ a=0,$ similar to the discussion for (\ref{eq:identity}), we have
\begin{equation}\label{eq:identity1}
\begin{aligned}
2\int_{r>r_{0},0<\theta<\pi}\eta|\nabla\omega|^{2}dxdy=&\int_{r>r_{0},0<\theta<\pi}\omega^{2}(\Delta\eta+\boldsymbol{w}\cdot\nabla\eta)dxdy\\&+\int_{r>r_{0},0<\theta<\pi}\nabla\eta\cdot (0,b)\omega^{2}dxdy.
\end{aligned}
\end{equation}
One verifies easily that there is a constant $ C $ independent of $ R $ such that
$$|\Delta\eta|\leq C,\quad |\nabla\eta|\leq\frac{C}{(\log r)^{\frac{1}{2}}}.$$
Noting that $ \eta=\frac{r}{\log r} $ for $ r_{1}<r<R, $ it follows  from (\ref{eq:identity1}) and (\ref{eq:velocitydecay})
\begin{equation*}
\begin{aligned}
&\int_{r_{1}<r<R,0<\theta<\pi}\frac{r}{\log r}|\nabla\omega|^{2}dxdy
\\&\leq \frac{1}{2}\int_{r>r_{0},0<\theta<\pi}\omega^{2}(\Delta\eta(r)+\boldsymbol{w}\cdot\nabla\eta(r))dxdy+\frac12\int_{r>r_{0},0<\theta<\pi}\nabla\eta(r)\cdot (0,b)\omega^{2}dxdy
\\&\leq C\int_{r>r_{0},0<\theta<\pi}\omega^{2}[1+\frac{|\boldsymbol{w}|}{(\log r)^{\frac{1}{2}}}]dxdy+C\int_{r>r_{0},0<\theta<\pi}\frac{1}{(\log r)^{\frac12}}\omega^{2}dxdy\\&\leq C(1+\frac{1}{(\log r_{0})^{\frac12}})\int_{r>r_{0},0<\theta<\pi}\omega^{2}dxdy
<\infty.
\end{aligned}
\end{equation*}
Letting $ R\to \infty, $ we obtain (\ref{eq:vor estimate}).

Using Lemma \ref{lem:vor estimate} one can improve the result of Lemma \ref{lem:vorticitydecay}.

{\bf  Proof of Theorem \ref{thm:vorticitydecay}.} Note that for $ 2^{n}>r_{0} $
\begin{equation*}
\begin{aligned}
&\int_{2^{n}}^{2^{n+1}}\frac{dr}{r}\int_{0}^{\pi}(r^{2}\omega^{2}+2\frac{r^{\frac{3}{2}}}{(\log r)^{\frac{1}{4}}}|\omega\omega_{\theta}|)d\theta
\\&\leq\int_{2^{n}<r<2^{n+1},0<\theta<\pi}(\omega^{2}+2\frac{r^{\frac{1}{2}}}{(\log r)^{\frac{1}{4}}}|\omega||\nabla\omega|)dxdy\\&\leq\int_{r>2^{n},0<\theta<\pi}(2\omega^{2}+2\frac{r}{(\log r)^{\frac{1}{2}}}|\nabla\omega)|^{2})dxdy.
\end{aligned}
\end{equation*}
Using (\ref{eq:vor estimate}) and proceeding exactly as in the proof of Lemma \ref{lem:vorticitydecay}, we obtain (\ref{eq:voricitydecay}). The proof is complete.

%

\noindent {\bf Acknowledgments.}
W. Wang was supported by NSFC under grant 12071054,  National Support Program for Young Top-Notch Talents and by Dalian High-level Talent Innovation Project (Grant 2020RD09).

\end{document}